\documentclass[twoside,11pt]{article}
\usepackage{amsmath,amssymb,amscd}

\newtheorem{lemma}{Lemma}
\newtheorem{prop}[lemma]{Proposition}
\newtheorem{theo}[lemma]{Theorem}
\newtheorem{cor}[lemma]{Corollary}
\newtheorem{defin}[lemma]{Definition}

\def\beginBeweis{{\it Proof . }}
\def\endBeweis{\hspace*{\fill}$\blacksquare$\\[2.5mm]}

\newcommand{\G}{\Gamma_g}
\newcommand{\bm}{\bigm|}

\begin{document}
\normalsize
\pagestyle{myheadings}
\markboth{R. Grieder}{$G$-Actions on Riemann Surfaces}
\noindent
{\Large\bf $G$-Actions on Riemann Surfaces and the\\
Associated Group of Singular Orbit Data}\\[3mm]
{\bf Ralph Grieder}\\[1mm]
{\small Department of Mathematics, Northwestern University, Evanston, 
IL 60208-2730}\\
{\small e-mail: ralph@math.nwu.edu }\\[8mm]

\begin{abstract}
Let $G$ be a finite group. To every smooth $G$-action on a
compact, connected and oriented Riemann surface we can associate its data of
singular orbits. The set of such data becomes an Abelian group
$\mathbb{B}_G$ under the $G$-equivariant connected sum. The map which
sends $G$ to $\mathbb{B}_G$ is functorial and carries many features of
the representation theory of finite groups. In this paper we will give
a complete computation of the group $\mathbb{B}_G$ for any finite
group $G$. 

There is a surjection from the $G$-equivariant cobordism group of
surface diffeomorphisms $\Omega_G$ to $\mathbb{B}_G$. We will prove
that the kernel of this surjection is isomorphic to
$H_2(G;\mathbb{Z})$. Thus $\Omega_G$ is an Abelian group extension of
$\mathbb{B}_G$ by $H_2(G;\mathbb{Z})$.

Finally we will prove that the group $\mathbb{B}_G$ contains only
elements of order two if and only if every complex character of $G$
has values in $\mathbb{R}$. This property shows a strong
relationship between the functor $\mathbb{B}$ and the
representation theory of finite groups.
\end{abstract}

\section{Introduction}\label{Intro}
Let $G$ be a finite group and $S_g$ a connected, oriented and compact Riemann
surface of genus $g$. 
One can assign to every diffeomorphic $G$-action on a surface $S_g$
its data of singular orbits and with this define the singular orbit data. 
In \cite{Gr3} the author defines a group structure on the set of all
singular orbit data and the resulting group is
denoted by $\mathbb{B}_G$ (see section \ref{sec:group} for definitions). 
In the same paper it is proven that the correspondence 
$$G\longmapsto\mathbb{B}_G$$
is a covariant functor. For finite subgroups $H$ of $G$ there are also
restriction and induction maps between the groups $\mathbb{B}_H$ and
$\mathbb{B}_G$, together with a double coset formula, thus the functor
$\mathbb{B}$ could be named {\it geometric representation
theory}.   

Another property of the group $\mathbb{B}_G$ which is proven in
\cite[section 6]{Gr3} states that there exists a surjection $\chi$
from $\Omega_G$, the $G$-equivariant cobordism group in dimension two,
onto $\mathbb{B}_G$ (see also section \ref{cobordism}).

The main purpose of this paper is to prove theorem \ref{isoB_G} which
gives a complete computation of the group $\mathbb{B}_G$ for any
finite group $G$.

In section \ref{cobordism} we will prove that the kernel of the above
map $\chi$ is isomorphic to $H_2(G;\mathbb{Z})$ 
and thus $\Omega_G$ is an Abelian group extension of $\mathbb{B}_G$ by
$H_2(G;\mathbb{Z})$. 

Moreover in section \ref{reptheory} we will show that the group
$\mathbb{B}_G$ contains only elements of order two if and only if
every complex character of $G$ has values in $\mathbb{R}$.
This property shows a strong relationship between the functor
$\mathbb{B}$ and the representation theory of finite groups.

Another reason for the interest in the group $\mathbb{B}_G$ comes from
the fact that it is possible to deduce information about the
cohomology of the mapping class group.

The mapping class group $\G$ is defined as follows. We write $\text{\it
Diffeo}_+(S_g)$ for the group of orientation preserving diffeomorphisms
of $S_g$ with the  $C^{\infty}$-topology, and $\text{\it Diffeo}^0_+(S_g)$
for the connected component of the identity. Then we have $\G=\text{\it
Diffeo}_+(S_g)/\text{\it Diffeo}^0_+(S_g)$. For a detailed discussion
of the mapping class group see \cite{Mi}.

The singular orbit data is also an invariant of a diffeomorphic
$G$-action up to isotopy, i.e., we can assign to every finite subgroup
$G$ of a mapping class group $\G$ its singular orbit data. 

The mapping class group acts on the first homology group of
$S_g$. This action preserves the intersection pairing and thus gives
rise to a symplectic representation $\eta:\G\rightarrow
Sp_{2g}(\mathbb{R})$. The unitary group $U(g)$ is a maximal compact
subgroup of $Sp_{2g}(\mathbb{R})$. Thus for any embedding
$\phi:G\rightarrow\G$ the map $\eta\circ\phi$ factors up to conjugation through
$U(g)\subset Sp_{2g}(\mathbb{R})$ and this unitary representation is denoted
by $\varphi(\phi)$. Hence $\varphi$ assigns to every embedding $\phi$ a
complex representation.

The important fact is now that the map $\varphi$ induces a group
homomorphism $\theta$ from $\mathbb{B}_G$ to a quotient of the complex
representation ring of $G$.
$$\theta:\mathbb{B}_G\rightarrow R_{\mathbb{C}}G/E_G$$
With this homomorphism $\theta$ we can now deduce results about the
cohomology of the mapping class group. 

In the case where $G$ is the cyclic group of prime
order $p$, it is shown in \cite{Gr1}
that the images of the symplectic classes $d_i\in H^{2i}
(BSp(\mathbb{R});\mathbb{Z})$ in $H^{2i}(\Gamma;\mathbb{Z})$ have
infinite order. Here we consider the stable situation, i.e., $\Gamma$
is the stable mapping class group. 

Another result in \cite{Gr2} states that one can embed polynomial
algebras in the cohomology of the stable mapping class group. For $p$
a regular prime, we have
$H^*(BSp(\mathbb{R});\mathbb{F}_p)\cong\mathbb{F}_p[d_1,d_2,\ldots]$, 
$deg\,d_i=2i$ and the map
$\eta^*:H^*(BSp(\mathbb{R});\mathbb{F}_p)\rightarrow H^*(\Gamma;\mathbb{F}_p)$ is injective on the polynomial algebra
$\mathbb{F}_p[(d_i)_{i\in J}]$, $J=\{i\in\mathbb{N}\,|\,i\equiv 1 \mod{2}\;
\text{or}\; i\equiv 0\mod{p-1}\}$.

\section{The Group of Singular Orbit Data}\label{sec:group}
Let $\phi$ be an embedding of any finite group $G$ into some mapping
class group $\G$. By Kerckhoff \cite{Ke} any such embedding can be
lifted to a homomorphism $G\rightarrow\text{\it Diffeo}_+(S_g)$, i.e.,
to an action of $G$ on the surface $S_g$ such that the elements of $G$
act by orientation preserving diffeomorphisms. Let $Gx$ denote the
orbit of $x\in S_g$ under the action of $G$. The orbit is called
singular if $|Gx|<|G|$, else
regular. If the orbit is singular, then there are elements of $G$
which stabilize the point $x\in S_g$ and $G_x$ denotes the
stabilizer of $G$ at $x$. It is proven by Accola \cite[Lemma 4.10]{Ac}
that the stabilizers are cyclic subgroups of $G$. Let $y$ be another element of
the orbit $Gx$. Thus there is an element $a\in G$ such that $ax=y$ and
the stabilizer of $y$ is conjugate to $G_x$, i.e., $aG_xa^{-1}=G_y$. 
As the elements of $G$ operate by orientation preserving
diffeomorphisms there are only finitely many singular orbits. Let
$x_i\in S_g$, $i=1,...,q$ be representatives of these orbits and $\nu_i$
the orders of the stabilizer groups $G_{x_i}$. Every group
$G_{x_i}$ has a generator $\gamma_i\in G$ such that $\gamma_i$ acts by
rotation through $2\pi/\nu_i$ on the tangent space at $x_i$. Similarly
$a\gamma_ia^{-1}$, $a\in G$, generates $aG_{x_i}a^{-1}$ and acts {\it
also} by rotation through $2\pi/\nu_i$ on the tangent space at $ax_i$.
Thus in order to collect information about the singular orbits, it is 
enough to consider the conjugacy classes of the elements $\gamma_i$.
Let $\hat{\gamma_i}$ denote the
conjugacy class of $\gamma_i$ in $G$. The {\bf singular orbit data}
of the embedding $\phi$ is then the unordered collection 
$$\bigl\{\hat{\gamma_1},\ldots,\hat{\gamma_q}\bigr\}_G.$$
By $W_G$ we will denote the set of all singular orbit data of $G$.
This data depends a priori on
the chosen lifting of $\phi$, but by \cite[Lemma 1]{Gr3} we know that
this data is well defined for an embedding of $G$. 
In the sequel
we will omit the subscript $_G$ if it is clear with respect
to which group the conjugacy classes are taken.

We have seen that every diffeomorphic $G$-action gives rise to a
singular orbit data. The question is now, which $q$-tuple of conjugacy
classes $\hat{\gamma_1},\ldots,\hat{\gamma_q}$ come form a $G$-action?
The next proposition answers this question.
\begin{prop}\label{sing}
An unordered $q$-tuple of conjugacy classes
$\hat{\gamma_1},\ldots,\hat{\gamma_q}$ is the singular orbit data of a
$G$-action if and only if $\gamma_1\cdots\gamma_q\in[G,G]$.\\
Here $[G,G]$ denotes the commutator subgroup of $G$.
\end{prop}
\beginBeweis
\cite[proposition 4]{Gr3}
\endBeweis
In the sequel we will denote by
$\bigl[\hat{\gamma_1},\ldots,\hat{\gamma_q}\bigr]_G$ an unordered
q-tuple of conjugacy classes of $G$ such that
$\gamma_1\cdots\gamma_q\in[G,G]$ and by $\Lambda_G$ the set of all
such q-tuples. We have now by proposition \ref{sing} a one to one
correspondence
$\bigl[\hat{\gamma_1},\ldots,\hat{\gamma_q}\bigr]_G\mapsto\bigl\{\hat{\gamma_1},\ldots,\hat{\gamma_q}\bigr\}_G$
between the elements of $\Lambda_G$ and $W_G$. In view of this
correspondence we will also call elements of $\Lambda_G$ singular
orbit data.  

Next we will define an addition, the {\bf $G$-equivariant connected sum}, on
the set of singular orbit data $W_G$.
Let $G$ act on a surface $S_g$ with singular orbit data
$\bigl\{\hat{\gamma_1},\ldots,\hat{\gamma_q}\bigr\}$ and on a surface
$S_h$ with singular orbit data
$\bigl\{\hat{\beta_1},\ldots,\hat{\beta_n}\bigr\}$. Find discs $D_1$
in $S_g$ and $D_2$ in $S_h$ such that $\bigl\{aD_j\bigr\}_{a\in G}$
are mutually disjoint for $j=1,2$. Then excise all discs
$\bigl\{aD_j\bigr\}_{a\in G}$, $j=1,2$, from $S_g$ and $S_h$ and take a
connected sum by matching $\partial(aD_1)$ to $\partial(aD_2)$ for all
$a\in G$. The resulting surface $S_{g+h+|G|-1}$ has $|G|$ tubes
joining $S_g$ and $S_h$. The actions of $G$ on $S_g$ and $S_h$ can be
extended to an action on $S_{g+h+|G|-1}$ by permuting the tubes. The
new action has a singular orbit data $\bigl\{\hat{\gamma_1},\ldots,
\hat{\gamma_q},\hat{\beta_1},\ldots,\hat{\beta_n}\bigr\}$.
This construction on surfaces defines an addition on the set
$W_G$.
$$\bigl\{\hat{\gamma_1},\ldots,\hat{\gamma_q}\bigr\}\oplus
\bigl\{\hat{\beta_1},\ldots,\hat{\beta_n}\bigr\}:=\bigl\{\hat{\gamma_1},
\ldots,\hat{\gamma_q},\hat{\beta_1},\ldots,\hat{\beta_n}
\bigr\}$$
After this geometric description we can give an algebraic
description of the addition on $\Lambda_G$.
Let $\bigl[\hat{\gamma_1},\ldots,\hat{\gamma_q}\bigr]_G$ and
$\bigl[\hat{\beta_1},\ldots,\hat{\beta_r}\bigr]_G$ be two elements of
$\Lambda_G$, the addition is now defined as follows:
$$\bigl[\hat{\gamma_1},\ldots,\hat{\gamma_q}\bigr]_G\oplus
\bigl[\hat{\beta_1},\ldots,\hat{\beta_r}\bigr]_G:=\bigl[\hat{\gamma_1},
\ldots,\hat{\gamma_q},\hat{\beta_1},\ldots,\hat{\beta_r}\bigr]_G.$$
With this addition we have only a commutative monoid structure on
$\Lambda_G$, respectively $W_G$, where the free 
actions represent the zero element. To obtain inverse elements we have to
introduce the following relations.

Suppose we have an action of $G$ on a surface $S_g$ with singular orbit data
$\bigl\{\hat{\gamma_1},\hat{\gamma_1}^{-1},$ $\hat{\gamma_2},\ldots,
\hat{\gamma_q}\bigr\}$. Let $\nu=|\bigl<\gamma_1\bigr>|$, then the
conjugacy class $\hat{\gamma_1}$ gives rise to a singular orbit with
representative $x$ such that $a\gamma_1a^{-1}$ acts by rotation through
$2\pi/\nu$ on the tangent space at $ax$, $a\in G$. On the other hand
$\hat{\gamma_1}^{-1}$ gives rise to another singular orbit with
representative $z$ such that $a\gamma_1a^{-1}$ acts by rotation through
$-2\pi/\nu$ on the tangent space at $az$, $a\in G$. Let $T$ be a set of
representatives for the $\bigl<\gamma_1\bigr>$ left cosets of $G$. Find discs
$D_1$ and $D_2$ around $x$ and $z$ respectively such that $D_j$ is fixed by
$\bigl<\gamma_1\bigr>$, $j=1,2$, and $\cup_{j=1,2}\cup_{t\in
T}\bigl\{tD_j\bigr\}$ are
mutually disjoint. Then excise all discs $\bigl\{tD_j\bigr\}_{t\in
T}$, $j=1,2$, from $S_g$ and connect the boundaries $\partial(tD_1)$
with $\partial(tD_2)$ by means of tubes $S^1\times[0,1]$ for every
$t\in T$. The resulting surface $S_{g+w}$ has $w=|G|/\nu$ new
handles. The action of $G$ on $S_g$ can be extended to $S_{g+w}$ by
permuting and rotating the new handles. This extended action yields
the singular orbit data
$\bigl\{\hat{\gamma_2},\ldots,\hat{\gamma_q}\bigr\}$. Pairs
of singular orbits which have opposite rotation on the tangent spaces
will be called cancelling pairs. The above process of eliminating
cancelling pairs will be called reduction and if there are no such
cancelling pairs left the singular orbit data is said to be in reduced form. 

Now we can define the relation.
$$\bigl\{\hat{\gamma_1},\ldots,\hat{\gamma_q}\bigr\}\sim
\bigl\{\hat{\beta_1},\ldots,\hat{\beta_n}\bigr\}:\Leftrightarrow
\left\{\begin{array}{l}
\text{The two singular orbit data}\\
\text{have the same reduced form}
\end{array}
\right.$$
This relation defines an equivalence relation on the set $W_G$ of
singular orbit data. 
$$\mathbb{W}_G:=W_G/\sim$$
The set $\mathbb{W}_G$ is not only a commutative monoid as $W_G$ but
contains also inverse elements and thus is a commutative group. The
inverse element of
$\bigl\{\hat{\gamma_1},\ldots,\hat{\gamma_q}\bigr\}$ is
$\bigl\{\hat{\gamma_1}^{-1},\ldots,\hat{\gamma_q}^{-1}\bigr\}$ and the
zero elements are cancelling pairs
$\bigl\{\hat{\gamma},\hat{\gamma}^{-1}\bigr\}$, 
$\gamma\in G$, and $\{\varnothing\}$ the free action.

We can also give a purely algebraic description of this group in terms
of $\Lambda_G$
\begin{equation}\label{B_G}
\mathbb{B}_G=\Lambda_G/\bigl<[\hat{\gamma},\hat{\gamma}^{-1}]_G
\,\bigm|\,\gamma\in G\bigr>
\end{equation}
with the inverse elements 
$\ominus[\hat{\gamma_1},\ldots,\hat{\gamma_q}]=[\hat{\gamma_1}^{-1},\ldots,
\hat{\gamma_q}^{-1}]$.
It is obvious from the definitions that the groups $\mathbb{B}_G$ and
$\mathbb{W}_G$ are canonically isomorphic (see also \cite[proposition
4]{Gr3}). In the sequel we will only use the 
notation $[\hat{\gamma_1},\ldots,\hat{\gamma_q}]$ for singular orbit
data and $\mathbb{B}_G$ for the {\bf group of singular orbit
data}. Even though elements of $\mathbb{B}_G$ consist of classes of
singular orbit data, we will by abuse of language also use the
notation $[\hat{\gamma_1},\ldots,\hat{\gamma_q}]$ for elements of
$\mathbb{B}_G$. But one has to keep in mind that with this notation
there is always a choice of representative involved. \\[1mm]
{\it Remark 1. }
All free $G$-actions represent the zero element in $\mathbb{B}_G$.
\\[1mm]
{\it Remark 2. }
The addition doesn't have any control on the genus. The genus can
become arbitrarily large.
\\[1mm]
Because of the purely algebraic description of $\mathbb{B}_G$ in
equation (\ref{B_G}) it is now possible to give a complete computation
of this group. This is done in theorem \ref{theoaisom} for finite Abelian
groups and for arbitrary finite groups in theorem \ref{isoB_G}.  

The map $\varphi$, which is defined in section \ref{Intro}, depends
only on the singular orbit data and the genus of the surface. Two
$G$-actions which have the same singular orbit data but not the same
genus are mapped under $\varphi$ to representations which differ only
by a rational representation. Moreover
the singular orbit data which form the relations,
$[\hat{\gamma},\hat{\gamma}^{-1}]_G$, $\gamma\in G$, are mapped under
$\varphi$ to rational representations and we obtain a well defined map
\begin{eqnarray}\label{eta}
\eta:\mathbb{B}_G&\rightarrow& R_{\mathbb{C}}(G)\\
  \alpha   & \mapsto
  &\varphi(\phi_{\alpha})-\overline{\varphi(\phi_{\alpha})}.\notag
\end{eqnarray}
Here $\phi_{\alpha}$ denotes a $G$-action which represents the
element $\alpha\in\mathbb{B}_G$ and $\overline{\varphi(\phi_{\alpha})}$
denotes the complex conjugate representation of
$\varphi(\phi_{\alpha})$. In addition one can prove that the
map $\eta$ is actually a group homomorphism.

The homomorphism $\eta\circ\chi$ is the $G$-signature defined by Atiyah and
Singer in \cite{AtSi}. This $G$-signature was used in the case
$G\cong\mathbb{Z}/p\mathbb{Z}$, $p$ a prime, by Ewing in \cite{Ew1}
and in the case $G\cong\mathbb{Z}/n\mathbb{Z}$, $n$ an integer, by
Edmonds/Ewing in \cite{EdEw} to prove their results.   

There is another way to define a $G$-signature with the help of the
map $\varphi$. By factoring out the image under $\varphi$ of all the
relations of 
$\mathbb{B}_G$, $\varphi$ induces a group homomorphism $\theta$.
$$\theta:\mathbb{B}_G\rightarrow R_{\mathbb{C}}(G)/E_G$$
Here $R_{\mathbb{C}}(G)/E_G$ denotes the complex representation ring
modulo the subgroup $E_G$. The map $\theta$ and the subgroup $E_G$ are
studied in some details in \cite[section 4]{Gr3}. In this paper we want
to concentrate on the group $\mathbb{B}_G$. However, in section
\ref{reptheory} we give a short overview of the map $\theta$ and the
subgroup $E_G$.

There are other interesting facts about the groups $\mathbb{B}_G$, which are
proven in \cite{Gr3}. E.g. for a homomorphism of finite groups
$f:H\rightarrow G$ there is a homomorphism
\begin{eqnarray*}
\mathbb{B}_f:\mathbb{B}_H & \rightarrow & \mathbb{B}_G \\
\bigl[\hat{\gamma_1},\ldots,\hat{\gamma_q}\bigr] & \mapsto &
\bigl[\widehat{f(\gamma_1)},\ldots,\widehat{f(\gamma_q)}\bigr]
\end{eqnarray*}
such that $\mathbb{B}_{f\circ g}=\mathbb{B}_f\circ\mathbb{B}_g$ and
$\mathbb{B}_{id}=id$
and thus the correspondence $\mathbb{B}:G\rightarrow\mathbb{B}_G$ is
functorial. This functor $\mathbb{B}$ carries many features of the
representation theory of finite groups. Let $H$ be a subgroup of $G$
and $i$ its inclusion. The inclusion induces the homomorphism $\mathbb{B}_i:
\mathbb{B}_H\rightarrow\mathbb{B}_G$ which is called
induction map. The terminology is motivated by the following
commutative diagram. 
\begin{eqnarray*}
\mathbb{B}_H & \stackrel{\mathbb{B}_i}{\rightarrow} & \mathbb{B}_G \\
\downarrow &  & \downarrow \\
R_{\mathbb{C}}H/E_H &
\stackrel{\widetilde{Ind}^G_H}{\rightarrow} & R_{\mathbb{C}}G/E_G
\end{eqnarray*}
Here $\widetilde{Ind}^G_H$ denotes the induction map restricted to the
quotient. 
Let $K$ be another subgroup of $G$, then there is another map
$\mathbb{B}res^G_K:\mathbb{B}_G\rightarrow\mathbb{B}_K$ called the
restriction map. It is defined by restricting the singular orbits of
$G$ to the subgroup $K$. This map satisfies also a commutative
diagram
\begin{eqnarray*}
\mathbb{B}_K & \stackrel{\mathbb{B}res^G_K}{\leftarrow} & \mathbb{B}_G \\
\downarrow & & \downarrow \\
R_{\mathbb{C}}K/E_K &
\stackrel{\widetilde{Res}^G_K}{\leftarrow} & R_{\mathbb{C}}G/E_G 
\end{eqnarray*}
where $\widetilde{Res}^G_K$ is the restriction map restricted to the
quotient. 
Furthermore the induction and restriction maps satisfy a double coset
formula (see \cite[proposition 5]{Gr3}). Thus the functor $\mathbb{B}$
describes a {\it geometric representation theory}.

Before we can turn to the next section we will give five examples to
illustrate the nature of the groups $\mathbb{B}_G$. We will also
prove two lemmata which will be useful throughout the paper. 

In the following examples the cyclic group of order $m$ will be
denoted by $C_m$.\\[1mm]
{\it Example 1. }
$G=C_2$; $\mathbb{B}_G\cong {0}$.\\[1mm]
{\it Example 2. }
Let $m$ be odd, $G=C_m=\bigl<x\bigr>$ then $\mathbb{B}_G\cong
\mathbb{Z}^{\frac{m-1}{2}}$ and a basis for $\mathbb{B}_G$ is given by
$\bigl[x,x^i,x^{m-i-1}\bigr]$, $i=1,\ldots,\frac{m-1}{2}$.\\[1mm]
{\it Example 3. }
Let $m$ be even, $G=C_m=\bigl<x\bigr>$ then $\mathbb{B}_G\cong
\mathbb{Z}^{\frac{m}{2}-1}$ and a basis for $\mathbb{B}_G$ is given by
$\bigl[x,x^i,x^{m-i-1}\bigr]$, $i=1,\ldots,\frac{m}{2}-1$.\\[1mm]
{\it Example 4. }
Let $p$ be an odd prime, $G=C_p\times C_p=\bigl<x\bigr>\times\bigl<y\bigr>$
then $\mathbb{B}_G\cong\mathbb{Z}^{\frac{p^2-1}{2}}$ and a basis for
$\mathbb{B}_G$ is given by $[x^j,y^i,x^{p-j}y^{p-i}]$,
$j=1,\ldots,$
 $(p-1)/2$, $i=1,\ldots,p-1$; $[x,x^k,x^{p-k-1}]$,
$k=1,\ldots,(p-1)/2$; $[y,y^l,y^{p-l-1}]$, $l=1,\ldots,(p-1)/2$.\\[2mm]
For $p$ the even prime, we have $\mathbb{B}_G\cong C_2$ and the
generator is $\bigl[x,y,xy\bigr]$.\\[1mm]
{\it Example 5. }
Let $G=S_3=C_3\rtimes C_2$, $C_3=\bigl<a\bigr>$ and
$C_2=\bigl<b\bigr>$. Then we have 
$[G,G]=C_3$, $\hat{a}=\{a,a^2\}$ and $\hat{b}=\{ab,a^2b,b\}$. The
singular orbit data are then $\bigl[\hat{a}\bigr]$,
$2\cdot\bigl[\hat{a}\bigr]=\bigl[\hat{a},\hat{a}\bigr]=\bigl[\hat{a},\hat{a^2}\bigr]=0$,
$\bigl[\hat{b},\hat{b}\bigr]=0$. Thus the group of singular orbit data
is generated by $\bigl[\hat{a}\bigr]$ and
$\mathbb{B}_{S_3}=\bigl<\bigl[\hat{a}\bigr]\bigr>\cong C_2$. 

As $\mathbb{B}_{C_2}\cong {0}$ the only maps which are of any interest are
the maps
$\mathbb{B}res^{S_3}_{C_3}:\mathbb{B}_{S_3}\rightarrow\mathbb{B}_{C_3}$
and $\mathbb{B}_i:\mathbb{B}_{C_3}\rightarrow\mathbb{B}_{S_3}$,
where $i$ is the inclusion $i:C_3\hookrightarrow S_3$. The maps are given by 
$\mathbb{B}res^{S_3}_{C_3}([\hat{a}]_{S_3})=[a,a^2]_{C_3}=0$ and 
$\mathbb{B}_i([a,a,a]_{C_3})=[\hat{a},\hat{a},\hat{a}]_{S_3}=
[\hat{a}]_{S_3}$. Thus $\mathbb{B}res^{S_3}_{C_3}$ is just the zero
map and $\mathbb{B}_i$ is the surjection of $\mathbb{Z}$ onto $C_2$.
\begin{lemma}\label{lemmatriple}
Any element of $\mathbb{B}_G$ can be written as a sum of triples 
$\bigl[\hat{x},\hat{y},\hat{z}\bigr]\in\mathbb{B}_G$.
\end{lemma}
\beginBeweis
Let $\bigl[\hat{x}_1,\ldots,\hat{x}_n\bigr]$ be any element of $\mathbb{B}_G$.
Then we can reduce the length by splitting off a triple.
\begin{align*}
\bigl[\hat{x}_1,\ldots,\hat{x}_n\bigr]=&\bigl[\hat{x}_1,\hat{x}_2,
\widehat{x_1x_2}^{-1}\bigr]\oplus\bigl[\widehat{x_1x_2},\hat{x}_3,\ldots,
\hat{x}_n\bigr]=\cdots\\
=&\bigoplus_{i=1}^{n-3}\bigl[\widehat{x_1\cdots x}_i,\hat{x}_{i+1},
\widehat{x_1\cdots x}_{i+1}^{-1}\bigr]\\
&\oplus\bigl[\widehat{x_1\cdots x}_{n-2},\hat{x}_{n-1},\hat{x}_n\bigr]
\end{align*}
\endBeweis
{\it Remark 3. }
The choice of the conjugacy class $\widehat{x_1\cdots x_i}$,
$i=2,\ldots,n-2$, in lemma \ref{lemmatriple}
is not unique. Every conjugacy class which is mapped to the same element in 
$G/[G,G]$ would also be a possible choice.
\begin{lemma}\label{lemmalinind}
For any group $G$ let $M=\{[\hat{x}_{i,1},\ldots,\hat{x}_{i,m_i}]\,|\,i=1,
\ldots,
m\}$, $Q=\{[\hat{z}_{i,1},\ldots,$ $\hat{z}_{i,q_i}]\,|\,i=1,
\ldots,
q\}$ and $N=\{[\hat{y}_{i,1},\ldots,\hat{y}_{i,n_i}]\,|\,i=1,\ldots,n\}$ be 
subsets of $\mathbb{B}_G$. Let $K$ be the subgroup generated by $M\cup
Q$ and $H$ the
subgroup generated by $N$. If the restrictions 
\begin{align}
\hat{x}_{i,m_i}&\neq\hat{x}_{i,m_i}^{-1}\;\;\text{for
all}\;\;i=1,\ldots,m \label{unique3}\\
\hat{x}_{i,m_i}&\neq\hat{x}_{k,l}^{\pm1}\;\;\text{for all}\;\;(i,m_i)\neq
(k,l)\;;\;i,k=1,\ldots,m\;;\;l=1,\ldots,m_k \label{unique1}\\
\hat{z}_{i,j}&=\hat{z}_{i,j}^{-1}\;\;\text{for
all}\;\;j=1,\ldots,q_i\;\;i=1,\ldots,q\label{unique4}\\
\hat{z}_{i,q_i}&\neq\hat{z}_{k,l}\;\;\text{for all}\;\;(i,q_i)\neq
(k,l)\;;\;i,k=1,\ldots,q\;;\;l=1,\ldots,q_k\label{unique5}\\
\hat{x}_{i,m_i}&\neq\hat{y}_{k,l}^{\pm1}\;\;\text{for all}\;\;i=1,\ldots,m\;;
\;k=1,\ldots,n\;;\;l=1,\ldots,n_k\label{unique2}\\
\hat{z}_{i,q_i}&\neq\hat{y}_{k,l}\;\;\text{for all}\;\;i=1,\ldots,q\;;
\;k=1,\ldots,n\;;\;l=1,\ldots,n_k\label{unique6}
\end{align}
apply, then the following hold.
\renewcommand{\theenumi}{\alph{enumi}}
\renewcommand{\labelenumi}{(\theenumi)}
\begin{enumerate}
\item There are no relations between the elements of $M\cup Q$ except
the obvious ones, $2\cdot\alpha=0\;,\;\forall\,\alpha\in Q$, thus
$K\cong\mathbb{Z}^m\oplus(\mathbb{Z}/2\mathbb{Z})^q$.  \label{indep}
\item The intersection of $K$ with $H$ is the 
trivial group and thus \\
$K\oplus H<\mathbb{B}_G$.\label{inters}
\end{enumerate}
\end{lemma}
\beginBeweis
Note that by equation (\ref{unique3}) and (\ref{unique4}) we have:
\begin{equation}
\hat{x}_{i,m_i}\neq\hat{z}_{k,j}\;\;\text{for
all}\;\;i=1,\ldots,m\;;\;j=1,\ldots,q_k\;;\;k=1,\ldots,q\label{unique7}
\end{equation}
and every element of $Q$ has order two.

First we prove statement (\ref{indep}). In a linear combination
\begin{equation}\label{lincomb1}
\bigoplus_{i=1}^{m}a_i\bigl[\hat{x}_{i,1},\ldots,\hat{x}_{i,m_i}\bigr]\oplus
\bigoplus_{i=1}^{q}\epsilon_i\bigl[\hat{z}_{i,1},\ldots,\hat{z}_{i,q_i}\bigr]
\;,\;a_i\in\mathbb{Z}\;,\;\epsilon_i=0,1
\end{equation}
the conjugacy classes $\hat{x}_{i,m_i}$, $i=1,\ldots,m$, cannot cancel because 
of equations (\ref{unique1}) and (\ref{unique7}) and thus they appear
exactly $a_i$ times. 
Consequently to have the linear combination (\ref{lincomb1}) 
equal zero, the coefficient $a_i$, $i=1,\ldots,m$, have to be
trivial. Furthermore by equation (\ref{unique5}) a similar argument
about the conjugacy classes $\hat{z}_{i,q_i}$, $i=1,\ldots,q$, shows that the
$\epsilon_i$, $i=1,\ldots,q$, have to be trivial. \\
Next we prove statement (\ref{inters}). Any nontrivial element $\eta$ of the 
intersection $K\cap H$ is 
a linear combination like in (\ref{lincomb1}). Thus at least one of the 
conjugacy 
classes $\hat{x}_{i,m_i}$, $i=1,\ldots,m$, or $\hat{z}_{i,q_i}$, $i=1,\ldots,q$ has to appear. On the other hand 
$\eta$ is also a linear combination of the $[\hat{y}_{i,1},\ldots,
\hat{y}_{i,n_i}]$, $i=1,\ldots,n$, and thus by equations
(\ref{unique2}) and (\ref{unique6}) the 
$\hat{x}_{i,m_i}$, $i=1,\ldots,m$ and $\hat{z}_{i,q_i}$,
$i=1,\ldots,q$ cannot appear which yields a contradiction. 
\endBeweis

\section{Finite Abelian Groups}\label{Abgrps}
In this section $G$ will always denote a finite Abelian group and the
maximal rank of an elementary Abelian $2$-subgroup of $G$ will be
denoted by $n_G$. First 
let $G$ be a cyclic group of prime power order $p^r$ with generator $x$,
then we can define the following subsets of $G$.
\begin{defin}\label{T_Gcyclic}
\begin{align*}
p=2\;\;,\; & T_G^+:=\{x^i\bm i=1,\ldots,p^r/2-1\}\\
      & S_G:=\{x^{p^r/2}\}\\
p\neq2\;\;,\;& T_G^+:=\{x^i\bm i=1,\ldots,(p^r-1)/2\}\\
        & S_G:=\emptyset
\end{align*}
\end{defin}
With these sets we can define $T_G^-:=\{x^{-1}\,|\,x\in T_G^+\}$ and thus
$T_G:=\{x\in G\,|\,x\neq x^{-1}\}$ is the disjoint union of $T_G^+$ with
$T_G^-$. Furthermore we have also $S_G=\{x\in
G\,|\,x=x^{-1}\;,\;x\neq1\}$. For a product  
$G=G_1\times G_2$ we introduce the following definitions.
\begin{defin}\label{T_Gproduct}
\begin{align}
T_G^+&:=T_{G_1}^+\cup T_{G_2}^+\cup \bigl(T_{G_1}^+\times S_{G_2}\bigr)\cup 
  \bigl(T_{G_1}^+\times T_{G_2}\bigr)\cup \bigl(S_{G_1}\times T_{G_2}^+\bigr)\notag\\
 &=T_{G_1}^+\cup T_{G_2}^+\cup \bigl(T_{G_1}^+\times G_2\setminus\{1\}\bigr)
   \cup\bigl(S_{G_1}\times T_{G_2}^+\bigr) \label{T_G^+}\\
S_G&:=S_{G_1}\cup S_{G_2}\cup \bigl(S_{G_1}\times S_{G_2}\bigr)\label{S_G}
\end{align}
In this definition elements $x$ of $G_1$ or $G_2$ are thought of as
elements $(x,1)$, respectively $(1,x)$, of $G$. 
\end{defin}
With the set $T_G^+$ we can again define the set $T_G^-:=\{x^{-1}\,|\,x\in
T_G^+\}$ and $T_G:=\{x\in G\,|\,x\neq x^{-1}\}$ is also the disjoint union
of $T_G^+$ with $T_G^-$. It is again true that $S_G=\{x\in
G\,|\,x=x^{-1}\;,\;x\neq1\}$. 

Every finite Abelian group $G$ is isomorphic to a product of cyclic
groups of prime power order. By definitions \ref{T_Gcyclic} and
\ref{T_Gproduct} we can define recursively for every such
factorization of $G$ together with a fixed choice of generator for
each factor, different sets $T_G^+$, $T_G^-$ and $S_G$. These different
sets however are isomorphic by the isomorphisms between the different
factorizations. In the sequel we will 
fix for every finite Abelian group one factorization into cyclic
subgroups of prime power order together with a generator for
each factor. Thus from now on $T_G^+$ and $S_G$ are fixed subsets of $G$. 
\begin{defin}\label{basiscyclic}
For $G$ a cyclic group of order $p^r$, $p$ a prime, $x$ the fixed generating 
element, we define
\begin{align*}
W_G&:=\Bigl\{\bigl[x,x^i,x^{p^r-1-i}\bigr]\in\mathbb{B}_G \bm \forall
i\in\mathbb{N}\; \mbox{such that}\;x^i\in
T_G^+\Bigr\}\\ 
V_G&:=\emptyset
\end{align*}
\end{defin}
\begin{prop}\label{cycbasis}
The set $W_G$ is a basis for $\mathbb{B}_G$ for any cyclic group $G$ of
prime power order.
\end{prop}
\beginBeweis
To prove that $W_G$ generates $\mathbb{B}_G$, it is enough by lemma
\ref{lemmatriple} to show that triples $\bigl[x^k,x^l,x^q\bigr]$ lie
in the span of $W_G$.

Let $\bigl[x^k,x^l,x^q\bigr]$ be any triple in $\mathbb{B}_G$ with
$k+l+q=sp^r$ and $1\leq k,l,q<p^r$, then we have the following
equations. (In the equations below the sums $q+k$, $q+t$ and $q+t+1$ will
always be modulo $p^r$.)
\begin{eqnarray*}
\bigl[x^k,x^l,x^q\bigr]&=&\bigl[x^k,x^l,x^{q+1},x^{p^r-1},x,x^{p^r-(q+1)},
x^q\bigr]\\
&=& \bigl[x^k,x^l,x^{q+1},x^{p^r-1}\bigr]\oplus\bigl[x,x^{p^r-(q+1)},
x^q\bigr]=\cdots= \\
&=&\bigl[x^k,x^l,x^{q+k},\underbrace{x^{p^r-1},\ldots,x^{p^r-1}}_k\bigr]
\oplus\bigoplus_{t=0}^{k-1}\bigl[x,x^{p^r-(q+t+1)},x^{q+t}\bigr] \\
&=&\ominus\bigl[\underbrace{x,\ldots,x}_k,x^{p^r-k}\bigr]\oplus
\bigoplus_{t=0}^{k-1}\bigl[x,x^{p^r-(q+t+1)},x^{q+t}\bigr]\\
&=&\bigoplus_{j=1}^{k-1}\ominus\bigl[x,x^j,x^{p^r-j-1}\bigr]
\bigoplus_{t=0}^{k-1}\bigl[x,x^{p^r-(q+t+1)},x^{q+t}\bigr]
\end{eqnarray*}
We will prove the linear independence only in the case $p\neq2$. Let 
$$\bigoplus_{i=1}^{(p^r-1)/2}a_i\bigl[x,x^i,x^{p^r-1-i}\bigr]\;\;,
\;\;a_i\in\mathbb{Z}$$ 
be any linear combination of elements in $W_G$. The elements $x^j$,
$j=2,\ldots,(p^r-3)/2$ and their inverse $x^{p^r-j}$ appear exactly
$a_j$, respectively $a_{j-1}$ times, in this linear combination, namely in 
$a_j\bigl[x,x^j,x^{p^r-1-j}\bigr]$ and in
$a_{j-1}\bigl[x,x^{j-1},x^{p^r-j}\bigr]$. The element $x^{(p^r-1)/2}$
appears exactly $2a_{(p^r-1)/2}$ times namely in
$a_{(p^r-1)/2}\bigl[x,x^{(p^r-1)/2},x^{(p^r-1)/2}\bigr]$ and its
inverse $a_{(p^r-3)/2}$ times in 
$a_{(p^r-3)/2}\bigl[x,x^{(p^r-3)/2},x^{(p^r+1)/2}\bigr]$. 

A necessary condition for the linear combination to be zero is such
that every element $x^j$, $j=1,\ldots,(p^r-1)/2$ cancels with its
inverse. This can only happen when
$$a_j=a_{j-1}\;\;,\;\;j=2,\ldots,(p^r-3)/2\;\;\mbox{and}\;\;a_{(p^r-3)/2}=2
a_{(p^r-1)/2}.$$
This however implies that the element $x$ appears $p^ra_{(p^r-1)/2}$
times but its 
inverse never. Thus the linear combination can only be zero when
$a_i=0$ for all $i=1,\ldots,(p^r-1)/2$.
\endBeweis
\begin{defin}\label{basisproduct}
Let $G$ be a product of two groups $G=G_1\times G_2$. An element $x\in G_i$, 
$i=1,2$, will also be considered as an element of $G$. Then we define 
\begin{align*} 
W_G:=&\Bigl\{\bigl[x,y,(x,y)^{-1}\bigr]\in\mathbb{B}_G \bm x\in T_{G_1}^+,y\in 
G_2\setminus\{1\}\;\mathrm{or}\;x\in S_{G_1},y\in T_{G_2}^+\Bigr\}\\
&\cup W_{G_1}\cup W_{G_2} \\
V_G:=&\Bigl\{\bigl[x,y,(x,y)\bigr]\in\mathbb{B}_G \bm x\in S_{G_1},y\in S_{G_2}
\Bigr\}\cup V_{G_1}\cup V_{G_2} 
\end{align*}
\end{defin}
{\it Remark 4. }
The set $V_G$ contains only elements of order two, i.e., for all 
$\bigl[x,y,(x,y)\bigr]\in V_G$ we have $2\cdot\bigl[x,y,(x,y)\bigr]=0$. On the 
other hand $W_G$ contains only elements of infinite order.
\\[1mm]
In the same way as we constructed the sets $T^+_G$ and $S_G$
recursively for any finite Abelian group $G$, we can construct by
definitions \ref{basiscyclic} and \ref{basisproduct} the subsets
$W_G$ and $V_G$ of $\mathbb{B}_G$ recursively for any finite Abelian group $G$.
\begin{theo}\label{theoabasis}
The set $W_G\cup V_G$ is a basis for $\mathbb{B}_G$ for any finite
Abelian group $G$.
\end{theo}
\beginBeweis
We will take advantage of the recursive definition of the sets $W_G$ and $V_G$
and prove the theorem by induction on the factorization of $G$ into
a product of groups.\\
For cyclic groups of prime power order the theorem follows from
proposition \ref{cycbasis}.\\
Let now $G$ be a product of $G_1$ and $G_2$ and assume that the theorem holds 
for the two groups $G_1$ and $G_2$. \\
First we prove that the set $W_G\cup V_G$ generates $\mathbb{B}_G$. By lemma 
\ref{lemmatriple} it is enough to prove that triples lie in the span of 
$W_G\cup V_G$. An arbitrary triple has the form $\bigl[(x_1,y_1),(x_2,y_2),
(x_1^{-1}x_2^{-1},y_1^{-1}y_2^{-1})\bigr]\in\mathbb{B}_G$, $x_i\in G_1$, 
$y_i\in G_2$, $i=1,2$, and we can write
\begin{multline*}
\bigl[(x_1,y_1),(x_2,y_2),(x_1^{-1}x_2^{-1},y_1^{-1}y_2^{-1})\bigr]=
\bigl[x_1^{-1},y_1^{-1},(x_1,y_1)\bigr]\\
\oplus\bigl[x_2^{-1},y_2^{-1},(x_2,y_2)\bigr]
\oplus\bigl[x_1x_2,y_1y_2,(x_1^{-1}x_2^{-1},y_1^{-1}y_2^{-1})\bigr]\\
\oplus\bigl[x_1,x_2,x_1^{-1}x_2^{-1}\bigr]
\oplus\bigl[y_1,y_2,y_1^{-1}y_2^{-1}\bigr].
\end{multline*}
The first three summand (resp. their inverse) of the right hand side of the 
equation are elements of $W_G\cup V_G$.
The forth summand is an element of $\mathbb{B}_{G_1}$ and thus by assumption 
a linear combination of elements in $W_{G_1}\cup V_{G_1}$, which is contained 
in $W_G\cup V_G$, the same holds for the fifth summand when we replace $G_1$ by
$G_2$. Note that whenever $x_i=1$ or $y_i=1$ then we ignore the elements
which are zero. \\
It remains to prove that there are no relations between the elements
of $W_G\cup V_G$ except the obvious ones, i.e., $2\alpha=0$,
$\forall\,\alpha\in V_G$.
To apply lemma \ref{lemmalinind} we have to specify what are the sets
$M$, $Q$ and $N$. In our situation we have that 
$M=W_G\setminus{(W_{G_1}\cup W_{G_2})}$, $Q=V_G\setminus{(V_{G_1}\cup
V_{G_2})}$ and $N=W_{G_1}\cup W_{G_2}\cup V_{G_1}\cup V_{G_2}$. 
By definition \ref{basisproduct} the equations (\ref{unique3}) to 
(\ref{unique6}) of lemma \ref{lemmalinind} are satisfied and therefore
there are no relations, except the obvious ones, between the elements
of $W_G\cup V_G\setminus{(W_{G_1}\cup W_{G_2}\cup V_{G_1}\cup
V_{G_2})}$ 
and also no relations with elements of 
$W_{G_1}\cup W_{G_2}\cup V_{G_1}\cup V_{G_2}$. 
On the other hand by assumption $W_{G_i}\cup V_{G_i}$ is a basis
of $\mathbb{B}_{G_i}$, $i=1,2$, and there are obviously no relations between 
$W_{G_1}\cup V_{G_1}$ and $W_{G_2}\cup V_{G_2}$.
\endBeweis
\begin{theo}\label{theoaisom}
$$\mathbb{B}_G \cong\mathbb{Z}^{|T_G^+|}\oplus(\mathbb{Z}/2\mathbb{Z})^{|S_G|-
n_G}$$
\end{theo}
\beginBeweis
By theorem \ref{theoabasis} and remark~4 we have 
$\mathbb{B}_G=\mathbb{Z}^{|W_G|}\oplus(\mathbb{Z}/2\mathbb{Z})^{|V_G|}$ thus it
suffices 
to prove $|W_G|=|T_G^+|$ and $|V_G|=|S_G|-n_G$. We will again proceed by
induction on 
the product structure.\\
For a cyclic group of prime power order $p^r$ we have by definitions
\ref{T_Gcyclic} and \ref{basiscyclic} 
\begin{align*}
p=2\;,\; & |W_G|=|T_G^+|\;, \\
 & |S_G|=1=n_G\; \;\text{and thus}\;\;|V_G|=0=|S_G|-n_G;
\end{align*}
\begin{align*}
p\neq 2\;,\; & |W_G|=|T_G^+|\;, \\
 & |S_G|=0=n_G\; \;\text{and thus}\;\;|V_G|=0=|S_G|-n_G.
\end{align*}
Let now $G$ be a product of $G_1$ and $G_2$ and assume that the theorem holds 
for the two groups $G_1$ and $G_2$, i.e.,  
\begin{xalignat*}{2}
|W_{G_1}|&=|T_{G_1}^+| &\qquad |V_{G_1}|&=|S_{G_1}|-n_{G_1} \\
|W_{G_2}|&=|T_{G_2}^+| &\qquad |V_{G_2}|&=|S_{G_2}|-n_{G_2}.
\end{xalignat*}
By definition \ref{basisproduct}, the equations (\ref{T_G^+}) and 
(\ref{S_G}) and the fact that $n_G=n_{G_1}+n_{G_2}$ we have
\begin{align*}
\begin{split}
|W_G| &=|T_{G_1}^+|\cdot(|G_2|-1)+ 
|S_{G_1}|\cdot|T_{G_2}^+|+|W_{G_1}|+|W_{G_2}|\\
   &=|T_{G_1}^+|\cdot(|G_2|-1)+ 
|S_{G_1}|\cdot|T_{G_2}^+|+|T_{G_1}^+|+|T_{G_2}^+| \\
  &=|T_G^+|,
\end{split}\\
\begin{split}
|V_G| &=|S_{G_1}|\cdot|S_{G_2}|+|V_{G_1}|+|V_{G_2}| \\
    &=|S_{G_1}|\cdot|S_{G_2}|+|S_{G_1}|-n_{G_1}+|S_{G_2}|-n_{G_2} \\
   & =|S_G|-n_G.
\end{split}
\end{align*}
\endBeweis
Later on in the Non-Abelian case it will be more convenient to replace 
$V_G$ by another 
basis. Let $D_2$ be the $2$-subgroup of $\mathbb{B}_G$, i.e., 
$D_2\cong(\mathbb{Z}/2\mathbb{Z})^{|S_G|-n_G}$. Furthermore let $H_2$ be the 
maximal elementary Abelian $2$-subgroup of $G$, i.e., 
$H_2\cong(\mathbb{Z}/2\mathbb{Z})^{n_G}$ and the elements 
$\{z_1,\ldots,z_{n_G}\}$ denote a generating set of $H_2$. By theorem 
\ref{theoaisom} and the fact that $S_{H_2}= H_2\setminus\{1\}$ 
and $T^+_{H_2}=\emptyset$ we have 
$$(\mathbb{Z}/2\mathbb{Z})^{|S_{H_2}|-n_G}=(\mathbb{Z}/2\mathbb{Z})^{|H_2|-1-
n_{H_2}}\cong\mathbb{B}_{H_2}<\mathbb{B}_G$$
\begin{prop}\label{H2}
$$\mathbb{B}_{H_2}= D_2<\mathbb{B}_G$$
\end{prop}
\beginBeweis
It suffices to prove that $H_2\setminus\{1\}= S_G$ .\\[2mm]
"$\subset$" Let $x\in H_2\setminus\{1\}\Rightarrow  x=x^{-1}\,,\,x\neq1\Rightarrow 
x\in S_G$\\
"$\supset$" Let $x\in S_G\setminus H_2$. The group generated by $H_2$ and $x$ 
has to be
an elementary Abelian $2$-group. As $H_2$ is maximal it follows $x\in H_2$
which contradicts the assumption.
\endBeweis
We have seen that $V_G$ is a basis for $D_2$ and thus by proposition \ref{H2} 
also a basis for $\mathbb{B}_{H_2}$. In the following we introduce another 
basis for $\mathbb{B}_{H_2}$.
\begin{defin}\label{L_G}
\setlength{\multlinegap}{0pt}
\begin{multline*}
L_G:=\Bigl\{\bigl[z_{i_1},\ldots,z_{i_s},z_{i_1}\cdots z_{i_s}\bigr]\in
\mathbb{B}_{H_2}\bigm| 
(i_1,\ldots,i_s)\;\text{an unordered } \\ \text{$s$-tuple, such that}\;
i_j=1,\ldots,n_G\;;\;i_j\neq i_k \;,\;\forall\,j\neq k\\
j,k=1,\ldots,s\;;\;2\leq s\leq n_G\Bigr\}
\end{multline*}
The elements $z_1,\ldots,z_{n_G}$ denote a generating set of $H_2$.
\end{defin}
\begin{prop}\label{LGbasis}
$L_G$ is a basis of $\mathbb{B}_{H_2}$.
\end{prop}
\beginBeweis
First we prove that $L_G$ spans $\mathbb{B}_{H_2}$. By lemma \ref{lemmatriple}
every element of $\mathbb{B}_{H_2}$ can be written as a sum of triples. 
The elements $z_i$, $i=1,\ldots,n_G$, generate $H_2$, thus an arbitrary triple 
has the form
$$\bigl[z_{i_1}\cdots z_{i_l},z_{j_1}\cdots z_{j_t}, z_{i_1}\cdots z_{i_l}\cdot
z_{j_1}\cdots z_{j_t}\bigr]$$
with $z_{i_k}\neq z_{i_h}$ and $z_{j_r}\neq z_{j_q}$ for all $k\neq h$
and $r\neq q$.
We can now write this element as a linear combination with elements from $L_G$.
\begin{multline*}
\bigl[z_{i_1}\cdots z_{i_l},z_{j_1}\cdots z_{j_t}, z_{i_1}\cdots z_{i_l}\cdot 
z_{j_1}\cdots z_{j_t}\bigr]\\
=\bigl[z_{i_1},\ldots,z_{i_l},z_{i_1}\cdots z_{i_l}\bigr]\oplus
\bigl[z_{j_1},\ldots,z_{j_t},z_{j_1}\cdots z_{j_t}\bigr]\\
\oplus\bigl[z_{i_1},\ldots,z_{i_l},z_{j_1},\ldots,z_{j_t},z_{i_1}\cdots 
z_{i_l}\cdot z_{j_1}\cdots z_{j_t}\bigr]
\end{multline*}
If $z_{i_k}=z_{j_r}$ for some $k$ and $r$ then reduce the last summand.\\
Next we prove that there are no relations. Each element $x\in H_2$ has
an unique 
presentation up to ordering $x=z_{i_1}\cdots z_{i_s}$,
$i_j=1,\ldots,n_G$, $i_j\neq i_k$, $\forall\,j\neq k$, 
$j,k=1,\ldots,s$, $1\leq s\leq n_G$. Thus there is a one to one
correspondence between 
the elements of $H_2$ which are the
product of at least two generators and the elements 
of $L_G$. Applying lemma \ref{lemmalinind}(\ref{indep}) with $Q=L_G$
and $M=\emptyset=N$ we deduce that there are no relations between the
elements of $L_G$ except the usual ones, i.e., $2\alpha=0$ for all
$\alpha\in L_G$.
\endBeweis
The conclusion is that both sets $L_G$ and $V_G$ are a basis for the same 
subgroup $\mathbb{B}_{H_2}$ of $\mathbb{B}_G$.
\section{Finite Groups}
Let $G$ be a finite group and $[G,G]$ denote the commutator subgroup of $G$. We
obtain a short exact sequence
$$1\rightarrow[G,G]\rightarrow G \stackrel{\phi}{\rightarrow}G'\rightarrow 1$$
where $G'=G/[G,G]$. We will write $\tilde{x}$ for the image of $x$ under 
$\phi$. 
The conjugacy class of $x$ in $G$ will be denoted by $\hat{x}$ and the set of 
conjugacy classes of $G$ by $\hat{G}$. By the notation $\hat{x}^{-1}$ we will 
mean $\widehat{x^{-1}}$. This makes sense since $\widehat{x^{-1}}=
\widehat{y^{-1}}$ if and 
only if $x$ and $y$ are conjugate. The homomorphism $\phi$ induces not only the
homomorphism $\mathbb{B}_{\phi}$ but also the well defined set map 
\begin{eqnarray*}
\hat{\phi}:\hat{G}&\rightarrow&\hat{G'}=G'.\\
\hat{x}&\mapsto&\tilde{x}
\end{eqnarray*}
In the sequel the symbol $\tilde{x}$ will denote an element of $G'$ and also 
the subset $x\cdot[G,G]\subset G$. Thus we can write for an element $\tilde{x}$
of $G'$
$$\tilde{x}=\bigcup_{i=0}^n\hat{x}_i\;,\;\hat{x}_i\in\hat{G}\;,\;i=0,\ldots,n.
$$
For every element $\tilde{x}\in G'$ we fix a numbering of its conjugacy classes
$\hat{x}_i$, $i=0,\ldots,n$, such that
\begin{itemize}
\item $\tilde{x}=\bigcup_{i=0}^n\hat{x}_i$,
\item if there are $k+1$ conjugacy classes with $\hat{x}_j=\hat{x}_j^{-1}$ then
$\hat{x}_i=\hat{x}_i^{-1}$ if and only if $i\in\{0,\ldots,k\}$,
\item if $\tilde{x}\neq\tilde{y}$ and $\tilde{x}^{-1}=\tilde{y}$ then
$\hat{x}_i^{-1}=\hat{y}_i$ for all $i=0,\ldots,n$.
\end{itemize}
{\it Remark 5. }
If $\tilde{x}\neq\tilde{x}^{-1}$ then $\hat{x}_i^{-1}\neq\hat{x}_j$,
$\forall\,j,i=0,\ldots,n$.
\\[1mm]
{\it Remark 6. }
If $\tilde{x}=\tilde{x}^{-1}$ and $\hat{x}_0^{-1}\neq\hat{x}_0$ then
for every $i=0,\ldots,n$, there exists some $j=0,\ldots,n$, $j\neq i$
with $\hat{x}_i^{-1}=\hat{x}_j$.
\\[1mm]
{\it Remark 7. }
If $\tilde{x}=\tilde{x}^{-1}$ and $\hat{x}_0^{-1}=\hat{x}_0$ then
$\hat{x}_i^{-1}=\hat{x}_i$, $\forall\,i=0,\ldots,k$, for some $0\leq
k\leq n$ and for every $k<i\leq n$, there exists some $k<j\leq n$, $j\neq i$
with $\hat{x}_i^{-1}=\hat{x}_j$.
\\[1mm]
With this numbering fixed we can define a well defined set map.
\begin{eqnarray*}
\hat{\psi}: G'&\rightarrow&\hat{G} \\
\tilde{x}&\mapsto&\hat{x}_0
\end{eqnarray*}
This map satisfies $\hat{\phi}\circ\hat{\psi}=id_{G'}$ and induces a map 
$\mathbb{B}_{\hat{\psi}}$ 
\begin{eqnarray*}
\mathbb{B}_{G'}&\stackrel{\mathbb{B}_{\hat{\psi}}}{\rightarrow} &\mathbb{B}_G\\
\bigl[\tilde{x},\tilde{y},\ldots\bigr]&\mapsto
&\bigl[\hat{\psi}(\tilde{x}),\hat{\psi}(\tilde{y}), 
\ldots\bigr].
\end{eqnarray*}
The map $\mathbb{B}_{\hat{\psi}}$ is not well defined, indeed if 
$\tilde{x}=\tilde{x}^{-1}$ but $\hat{x}_0\neq\hat{x}_0^{-1}$ then
$$0=\mathbb{B}_{\hat{\psi}}(0)=\mathbb{B}_{\hat{\psi}}(\bigl[\tilde{x},
\tilde{x}\bigr])=\bigl[\hat{x}_0,\hat{x}_0\bigr]\neq0.$$
We collect the images of the elements where the map $\mathbb{B}_{\hat{\psi}}$ 
fails to be well defined: 
$P_G=\left\{\bigl[\hat{x}_0,\hat{x}_0\bigr]\bigm| 
\tilde{x}=\tilde{x}^{-1}\;,\;\hat{x}_0\neq\hat{x}_0^{-1}\right\}$.
\begin{prop}\label{P_G}
The map $\mathbb{B}_{\hat{\psi}}$ is well defined and linear up to elements of
$P_G$.
\end{prop}
\beginBeweis
The only relations in $\mathbb{B}_{G'}$ which are not satisfied in
$\mathbb{B}_G$, under the map $\mathbb{B}_{\hat{\psi}}$, are given by 
$\left\{\bigl[\tilde{x},\tilde{x}\bigr]\bigm| 
\tilde{x}=\tilde{x}^{-1}\;,\;\hat{x}_0\neq\hat{x}_0^{-1}\right\}$. 
Thus up to elements of $P_G$ the map $\mathbb{B}_{\hat{\psi}}$ is well defined.

Let $\bigl[\,_1\tilde{x},\ldots,\,_l\tilde{x}\bigr]$ and $\bigl[\,_1\tilde{y},
\ldots,\,_k\tilde{y}\bigr]$ be two elements in $\mathbb{B}_{G'}$. (We
introduce here the indexing on the left, because we will need also the
indexing on the right later on.) Their sum is
given by $\bigl[\,_1\tilde{x},\ldots,\,_l\tilde{x},\,_1\tilde{y},
\ldots,\,_k\tilde{y}\bigr]$ up to cancelling pairs. Under the map 
$\mathbb{B}_{\hat{\psi}}$ we obtain
$$\bigl[\,_1\hat{x}_0,\ldots,\,_l\hat{x}_0\bigr]\oplus\bigl[\,_1\hat{y}_0,
\ldots,\,_k\hat{y}_0\bigr]\equiv\bigl[\,_1\hat{x}_0,\ldots,\,_l\hat{x}_0,
\,_1\hat{y}_0,\ldots,\,_k\hat{y}_0\bigr]$$
up to the image of cancelling pairs in $\mathbb{B}_{G'}$. The only 
pairs which cancel in $\mathbb{B}_{G'}$ but their images under
$\mathbb{B}_{\hat{\psi}}$ do not cancel in $\mathbb{B}_G$ are the 
elements $\{\bigl[\tilde{x},\tilde{x}\bigr]\bigm| 
\tilde{x}=\tilde{x}^{-1}\;,$ $\hat{x}_0\neq\hat{x}_0^{-1}\}$. 
Thus up to elements of $P_G$ the map $\mathbb{B}_{\hat{\psi}}$ is linear.
\endBeweis
Note that whenever $\bigl[\hat{x}_0,\hat{y}_0,\ldots\bigr]$ is an element of 
$\mathbb{B}_G$ and $\hat{\phi}(\hat{x}_0)=\hat{\phi}(\hat{x}_i)$ and
$\hat{\phi}(\hat{y}_0)=\hat{\phi}(\hat{y}_i)$ then
$\bigl[\hat{x}_i,\hat{y}_j,\ldots\bigr]$ is also an element
of $\mathbb{B}_G$ for any $i$ and $j$. 
Moreover they have the same image under $\mathbb{B}_{\phi}$; 
$\mathbb{B}_{\phi}(\bigl[\hat{x}_0,\hat{y}_0,\ldots\bigr])=\mathbb{B}_{\phi}
(\bigl[\hat{x}_i,\hat{y}_j,\ldots\bigr])=\bigl[\tilde{x},\tilde{y},\ldots
\bigr]$. \\
We define again the sets $S_G$, $T_G$ and 
$T_G^{\pm}$ but with more conditions.
\begin{align*}
S_G&=\bigl\{\hat{x}\in\hat{G}\bigm|\hat{x}^{-1}=\hat{x}\;,\;x\neq1\bigr\}\\
T_G&=\bigl\{\hat{x}\in\hat{G}\bigm|\hat{x}^{-1}\neq\hat{x}\bigr\}
\end{align*}
The set $T_G$ is again the disjoint union of $T_G^+$ with $T_G^-$ satisfying 
the following conditions.
\begin{enumerate}
\item If $\hat{x}\in T_G^+$ then $\hat{x}^{-1}\in T_G^-$.
\item If $\hat{x}\in T_G^+$ with $\tilde{x}\neq\tilde{x}^{-1}$ then $\hat{y}\in
T_G^+$ for every $y\in x\cdot[G,G]$.
\item If $\tilde{x}=\tilde{x}^{-1}$ with $\hat{x}_0\neq\hat{x}_0^{-1}$
then $\hat{x}_0\in T_G^+$.
\end{enumerate}
Note that if $\hat{x}_i\in S_G$ for some $i$ then $\tilde{x}=\tilde{x}^{-1}$ 
and by the way we fixed the numbering of the conjugacy classes it follows
$\hat{x}_0\in S_G$. \\
With $H'_2$ we denote the maximal elementary Abelian 2-subgroup of
$G'$ and with  
$n_{G'}$ its rank, i.e., $H'_2\cong(\mathbb{Z}/2\mathbb{Z})^{n_{G'}}$. Let 
$$M_{G'}=\left\{\tilde{x}\in G'\;|\; \hat{x}_0^{-1}=\hat{x}_0\right\}$$
be a subset of $H'_2$ and $K'$ the subgroup generated by $M_{G'}$. For the rank
of $K'$ we write $n_{K'}$ and then we have
$$(\mathbb{Z}/2\mathbb{Z})^{n_{K'}}\cong K'<H'_2\cong (\mathbb{Z}/2\mathbb{Z})
^{n_{G'}}.$$
Now we have all the ingredients to construct a basis of $\mathbb{B}_{G}$ and 
with this basis to prove theorem \ref{isoB_G}.
\renewcommand{\theenumi}{\roman{enumi}}
\renewcommand{\labelenumi}{(\theenumi)}
\begin{align*}
N_1&=\Bigl\{\bigl[\hat{x}\bigr]\bigm|\tilde{x}=1\;,\;\hat{x}\in T_G^+\cup
  S_G\Bigr\}\\
N_2&=\Bigl\{\bigl[\hat{x}_0,\hat{x}_i^{-1}\bigr]\bigm|\tilde{x}\neq
  \tilde{x}^{-1}\;,\;\hat{x}_0\in T_G^+\;,\;1\leq i\leq n\Bigr\}\\
N_3&=\Bigl\{\bigl[\hat{x}_0,\hat{x}_i\bigr]\bigm|\tilde{x}=\tilde{x}^{-1}\;,\;
  \tilde{x}\neq1\;,\;\hat{x}_0\in T_G^+\;,\;1\leq i\leq
  n\;\,\mbox{such that}\;\hat{x}_i\in T_G^+\Bigr\}\\ 
N_4&=\Bigl\{\bigl[\hat{x}_0,\hat{x}_i\bigr]\bigm|\tilde{x}=\tilde{x}^{-1}\;,\; 
  \tilde{x}\neq1\;,\;\hat{x}_0\in S_G\;,\;1\leq i\leq n\;\,\mbox{such that}\; 
  \hat{x}_i\in T_G^+\cup S_G\Bigr\}\\
N_5&=\Bigl\{\bigl[\hat{x}_0,\hat{y}_0,\hat{z}_0\bigr]=
  \mathbb{B}_{\hat{\psi}}(\bigl[\tilde{x},\tilde{y},\tilde{z}\bigr])
  \in\mathbb{B}_{G}\bigm|\bigl[\tilde{x},\tilde{y},\tilde{z}\bigr]\in 
  W_{G'}\Bigr\}
\end{align*}
For the last sets $N_6$ and $N_7$ we need the basis $L_{G'}$ of
$\mathbb{B}_{H'_2}$ with a special choice of the generating elements
of $H'_2$. Choose a generating set $\{\,_1\tilde{z},\ldots,\,_{n_{K'}}
\tilde{z}\}$ of $K'$ such that $\,_i\tilde{z}\in M_{G'}$,
$i=1,\ldots,n_{K'}$. Now find another set $\{\,_{(n_{K'}+1)}\tilde{z},
\ldots,\,_{n_{H'_2}}\tilde{z}\}$
such that the union of both $\{\,_1\tilde{z},\ldots,\,_{n_{H'_2}}\tilde{z}\}$ 
generates the group $H'_2$. With this generating set and definition \ref{L_G}
we can construct $L_{G'}$:
\setlength{\multlinegap}{0pt}
\begin{multline*}
L_{G'}=\Bigl\{\bigl[\,_{i_1}\tilde{z},\ldots,\,_{i_s}\tilde{z},
\,_{i_1}\tilde{z}\cdots\,_{i_s}\tilde{z}\bigr]\in
\mathbb{B}_{H'_2}\bigm| 
(i_1,\ldots,i_s)\;\text{an unordered } \\ \text{$s$-tuple, such that}\;
i_j=1,\ldots,n_{G'}\;;\;i_j\neq i_k \;,\;\forall\,j\neq k\\
\forall\, j,k=1,\ldots,s\;;\;2\leq s\leq n_{G'}\Bigr\}.
\end{multline*}
Now we are ready to define the last two sets which will be part of a basis
of~$\mathbb{B}_{G}$.
\begin{equation*}
\begin{split}
N_6=&\Bigl\{\bigl[\,_{i_1}\hat{z}_0,\ldots,\,_{i_s}\hat{z}_0,
\,_{i_1}\hat{z}_0\cdots\,_{i_s}\hat{z}_0\bigr]=
\mathbb{B}_{\hat{\psi}}(\bigl[\,_{i_1}\tilde{z},\ldots,\,_{i_s}\tilde{z},
\,_{i_1}\tilde{z}\cdots\,_{i_s}\tilde{z}\bigr])\in\mathbb{B}_{G}\bigm|\\ 
& \;\;\;\bigl[\,_{i_1}\tilde{z},\ldots,\,_{i_s}\tilde{z},\,_{i_1}\tilde{z}\cdots\,_{i_s}\tilde{z}\bigr]\in L_{G'}\Bigr\} \\
N_7=&\Bigl\{\bigl[\,_i\hat{z}_0,\,_i\hat{z}_0\bigr]=
\bigl[\hat{\psi}(\,_i\hat{z}),\hat{\psi}(\,_i\hat{z})\bigr]\in
\mathbb{B}_{G}\bigm|i=n_{K'}+1,\ldots,n_{G'}\Bigr\}\\
\end{split}
\end{equation*}
Note that the product of conjugacy classes is not defined. The symbol 
$_{i_1}\hat{z}_0\cdots\,_{i_s}\hat{z}_0$ stands for the image of 
$_{i_1}\tilde{z}\cdots\,_{i_s}\tilde{z}$ under the map $\hat{\psi}$.
We haven't added the elements $\bigl[\,_i\hat{z}_0, _i\hat{z}_0\bigr]$, 
$i=1,\ldots,n_{K'}$, as they are all trivial. Indeed the elements $_i\tilde{z}$, 
$i=1,\ldots,n_{K'}$, belong to $M_{G'}$ and thus $_i\hat{z}_0= 
\,_i\hat{z}_0^{-1}$.
\begin{prop}\label{pair}
The elements $\bigl[\hat{x}_0,\hat{x}_0\bigr]\in P_G$ lie all in the
subgroup generated by $N_6\cup N_7$. 
\end{prop}
\beginBeweis
Let $\bigl[\hat{x}_0,\hat{x}_0\bigr]$ be an element of $P_G$, then
$\hat{\phi}(\hat{x}_0)=\tilde{x}$ is an element of $H_2'$ and thus
$\tilde{x}=\,_{i_1}\tilde{z}\cdots\,_{i_t}\tilde{z}$ for some
$i_j=1,\ldots,n_{H'_2}$, $j=1,\ldots,t$, furthermore
$\hat{\psi}(\tilde{x})=\hat{x}_0$ Now we have two cases either
$t=1$ or $t>1$.

For $t=1$ we have $_{i_1}\hat{z}_0=\hat{\psi}(\,_{i_1}\tilde{z})=
\hat{\psi}(\hat{\phi}(\hat{x}_0))=\hat{x}_0$, thus $n_{K'}<i_1\leq
n_{H'_2}$ and by the definition of $N_7$ it follows
$\bigl[\hat{x}_0,\hat{x}_0\bigr]\in N_7$.

For $t>1$ we obtain the following equation;
\begin{eqnarray*}
\bigl[\hat{x}_0,\hat{x}_0\bigr]&=&2\cdot\bigl[\,_{i_1}\hat{z}_0,\ldots,
\,_{i_t}\hat{z}_0,\hat{x}_0\bigr]\bigoplus^t_{j=1}\ominus\bigl[\,_{i_j}\hat{z}_0,
\,_{i_j}\hat{z}_0\bigr]\\
&=&2\cdot\bigl[\,_{i_1}\hat{z}_0,\ldots,\,_{i_t}\hat{z}_0,\,
_{i_1}\hat{z}_0\cdots\,_{i_t}\hat{z}_0\bigr]\bigoplus^t_{j=1}\ominus
\bigl[\,_{i_j}\hat{z}_0,\,_{i_j}\hat{z}_0\bigr]
\end{eqnarray*}
First note that the elements $\bigl[\,_{i_j}\hat{z}_0,\,_{i_j}\hat{z}_0\bigr]$ 
are zero whenever $1\leq i_j\leq n_{K'}$ and otherwise they belong
to $N_7$. On the other hand the element $\bigl[\,_{i_1}\hat{z}_0, 
\ldots,\,_{i_t}\hat{z}_0,\,_{i_1}\hat{z}_0\cdots\,_{i_t}\hat{z}_0\bigr]$
belongs to $N_6$ and the proposition is proven.
\endBeweis
\begin{prop}\label{spanB_G}
The set $\cup_{i=1}^7N_i$ generates $\mathbb{B}_G$.
\end{prop}
\beginBeweis
By lemma \ref{lemmatriple} it suffices to prove that any triple
$\bigl[\hat{x},\hat{y},\hat{z}\bigr]$ is a linear combination of elements in
$\cup_{i=1}^7N_i$. 

We fix the following notation: $\hat{\phi}(\hat{x})=\tilde{x}$,
$\hat{\phi}(\hat{y})=\tilde{y}$, $\hat{\phi}(\hat{z})=\tilde{z}$ and
$\hat{\psi}(\tilde{x})=\hat{x}_0$, $\hat{\psi}(\tilde{y})=\hat{y}_0$,
$\hat{\psi}(\tilde{z})=\hat{z}_0$. 

We have to distinguish four cases:
\begin{enumerate}
\item
Let $x,y,z\in[G,G]$ then we can write $\bigl[\hat{x},\hat{y},\hat{z}\bigr]= 
\bigl[\hat{x}\bigr]\oplus\bigl[\hat{y}\bigr]\oplus\bigl[\hat{z}\bigr]$
and thus it is generated by $N_1$.
\item
Let $x,y\in[G,G]$ then we have $x\cdot y\cdot z\in[G,G]$ and thus
$z\in[G,G]$ and we reduced this case to the first one.
\item
Let $x\in[G,G]$ and $y,z\notin[G,G]$ then we can write
$\bigl[\hat{x},\hat{y},\hat{z}\bigr]= 
\bigl[\hat{x}\bigr]\oplus\bigl[\hat{y},\hat{z}\bigr]$ where
$\bigl[\hat{x}\bigr]$ lies in the span of $N_1$. Thus it is enough to
show that $\bigl[\hat{y},\hat{z}\bigr]$ is in the span of
$\cup_{i=1}^7N_i$. Note that in this case $\tilde{z}=\tilde{y}^{-1}$.
Now assume that $\bigl[\hat{y},\hat{z}\bigr]$ isn't a cancelling pair,
then we have again different cases: 
\begin{enumerate}
\item
Let $\tilde{z}\neq\tilde{z}^{-1}$ which implies $\hat{z}_0\in T_G$ and
$\hat{z}_0=\hat{y}_0^{-1}$, then we have either $\hat{z},\hat{z}_0\in
T^+_G$ or $\hat{y},\hat{y}_0\in T^+_G$. 
\begin{itemize}
\item
If $\hat{z},\hat{z}_0\in T^+_G$ then of course $\hat{y},\hat{y}_0\in
T^-_G$ and we obtain
$\bigl[\hat{y},\hat{z}\bigr]=\ominus\bigl[\hat{z}_0,\hat{z}^{-1}\bigr]\oplus
\bigl[\hat{z}_0,\hat{y}\bigr]$. The last two summands are elements of $N_2$.
\item
If $\hat{y},\hat{y}_0\in T^+_G$ then equivalently $\hat{z},\hat{z}_0\in T^-_G$ and we get
$\bigl[\hat{y},\hat{z}\bigr]=\ominus\bigl[\hat{y}_0,\hat{y}^{-1}\bigr]\oplus
\bigl[\hat{y}_0,\hat{z}\bigr]$. The last two summands are again
elements of $N_2$.
\end{itemize}
Note that whenever $\hat{z}_0=\hat{z}$ or $\hat{y}_0=\hat{y}$ then we
ignore the elements which are zero.
\item
Let $\tilde{z}=\tilde{z}^{-1}$, which implies $\hat{z}_0=\hat{y}_0$,
and assume that $\hat{z}_0\in T^+_G$. These assumptions lead to the
following cases.
\begin{itemize}
\item
$\hat{y},\hat{z}\in T^+_G$ $\Rightarrow$ $\bigl[\hat{y},\hat{z}\bigr]=
\bigl[\hat{y}_0,\hat{y}\bigr]\oplus\bigl[\hat{z}_0,\hat{z}\bigr]\ominus
\bigl[\hat{y}_0,\hat{z}_0\bigr]$
\item
$\hat{y}\in T^+_G$, $\hat{z}\in T^-_G$ $\Rightarrow$ $\bigl[\hat{y},
\hat{z}\bigr]=\bigl[\hat{y}_0,\hat{y}\bigr]\ominus
\bigl[\hat{z}_0,\hat{z}^{-1}\bigr]$
\item
$\hat{y},\hat{z}\in T^-_G$ $\Rightarrow$ $\bigl[\hat{y},\hat{z}\bigr]=
\ominus\bigl[\hat{y}_0,\hat{y}^{-1}\bigr]\ominus
\bigl[\hat{z}_0,\hat{z}^{-1}\bigr]\oplus\bigl[\hat{y}_0,\hat{z}_0\bigr]$
\end{itemize}
We see that in the above three cases the element
$\bigl[\hat{y},\hat{z}\bigr]$ is generated by $N_3$ and by proposition
\ref{pair} also by the set $N_6\cup N_7$.
\item
Let $\tilde{z}=\tilde{z}^{-1}$, which implies $\hat{z}_0=\hat{y}_0$,
and assume that $\hat{z}_0\in S_G$. These assumptions lead to the
following cases.
\begin{itemize}
\item
$\hat{y},\hat{z}\in T^+_G\cup S_G$ $\Rightarrow$ $\bigl[\hat{y},\hat{z}\bigr]=
\bigl[\hat{y}_0,\hat{y}\bigr]\oplus\bigl[\hat{z}_0,\hat{z}\bigr]$
\item
$\hat{y}\in T^+_G\cup S_G$, $\hat{z}\in T^-_G$ $\Rightarrow$ $\bigl[\hat{y},\hat{z}\bigr]=
\bigl[\hat{y}_0,\hat{y}\bigr]\ominus\bigl[\hat{z}_0,\hat{z}^{-1}\bigr]$
\item
$\hat{y},\hat{z}\in T^-_G$ $\Rightarrow$ $\bigl[\hat{y},\hat{z}\bigr]=
\ominus\bigl[\hat{y}_0,\hat{y}^{-1}\bigr]\ominus\bigl[\hat{z}_0,\hat{z}^{-1}\bigr]$
\end{itemize}
For these three cases we can deduce that the element 
$\bigl[\hat{y},\hat{z}\bigr]$ is generated by $N_4$. Note that
whenever $\hat{z}_0=\hat{z}$ or $\hat{y}_0=\hat{y}$ then we 
ignore the elements which are zero.
\end{enumerate}
\item
For the last case we have $x,y,z\notin[G,G]$ which implies:
$$\bigl[\hat{x},\hat{y},\hat{z}\bigr]=\bigl[\hat{x}_0,\hat{y}_0,\hat{z}_0\bigr]
\oplus\bigl[\hat{x}^{-1}_0,\hat{x}\bigr]\bigl[\hat{y}^{-1}_0,\hat{y}\bigr]\bigl[\hat{z}^{-1}_0,\hat{z}\bigr].$$
The last three summands are pairs and by the third case lie in the
span of $\cup_{i=1}^7N_i$ and so we have to consider the element
$\bigl[\hat{x}_0,\hat{y}_0,\hat{z}_0\bigr]$.
\begin{eqnarray*}
\bigl[\hat{x}_0,\hat{y}_0,\hat{z}_0\bigr]&=&\mathbb{B}_{\hat{\psi}}\circ
\mathbb{B}_{\phi}(\bigl[\hat{x}_0,\hat{y}_0,\hat{z}_0\bigr])\\
&=&\mathbb{B}_{\hat{\psi}}(\bigl[\tilde{x},\tilde{y},\tilde{z}\bigr])
\end{eqnarray*}
The element $\bigl[\tilde{x},\tilde{y},\tilde{z}\bigr]$ belongs to
$\mathbb{B}_{G'}$ and thus by theorem \ref{theoabasis} and proposition
\ref{LGbasis} is a linear combination of elements
$\{\theta_i\}_{i\in I}$ in $W_{G'}\cup L_{G'}$. By proposition
\ref{P_G} the map $\mathbb{B}_{\hat{\psi}}$ is linear up to elements
of $P_G$ and so we get:
\begin{eqnarray*}
\bigl[\hat{x}_0,\hat{y}_0,\hat{z}_0\bigr]&=&\mathbb{B}_{\hat{\psi}}(\bigl[\tilde{x},\tilde{y},\tilde{z}\bigr])\\
&=&\bigoplus_{i\in I}\mathbb{B}_{\hat{\psi}}(\theta_i)\bigoplus_l\lambda_l
\end{eqnarray*}
where $\lambda_l$ are elements of $P_G$. By proposition \ref{pair}
the elements of $P_G$ are generated by $N_6\cup N_7$ and the elements 
$\mathbb{B}_{\hat{\psi}}(\theta_i)$ lie in $N_5$ or $N_6$.
\end{enumerate}
We have now proven the proposition by showing that any triple
$\bigl[\hat{x},\hat{y},\hat{z}\bigr]$ is generated by elements in
$\cup_{i=1}^7N_i$. 
\endBeweis
\begin{prop}\label{linunabB_G}
There are no relations between the elements of $\cup_{i=1}^7N_i$
except for some elements in $N_1$, $N_4$ and $N_6$ which have order two.
\end{prop}
\beginBeweis
First we apply lemma \ref{lemmalinind} to the sets $M\cup Q=N_1\cup N_2\cup N_3\cup N_4$
and $N=N_5\cup N_6\cup N_7$ to deduce that there are no relations
between the elements of $N_1\cup N_2\cup N_3\cup N_4$ and also no
relations with elements of $N_5\cup N_6\cup N_7$, besides of course the
two torsion. \\
Next we consider the set $N_5$. Suppose we could express an element
$\alpha$ of the
group generated by $N_6\cup N_7$ as a linear combination 
$$\alpha=\bigoplus_{i=1}^ka_i\bigl[\,_i\hat{x}_0,\,_i\hat{y}_0,\,_i\hat{z}_0
\bigr]\;\;,\;\;a_i\in\mathbb{Z}$$
where $\bigl[\,_i\hat{x}_0,\,_i\hat{y}_0,\,_i\hat{z}_0\bigr]$,
$i=1,\ldots,k$, are $k$ different elements of $N_5$.
The map $\mathbb{B}_{\phi}$ sends $\alpha$ to an element of the
group generated by $L_{G'}$ and the linear combination to
$$\mathbb{B}_{\phi}(\alpha)=\bigoplus_{i=1}^ka_i\bigl[\,_i\tilde{x},
\,_i\tilde{y},\,_i\tilde{z}\bigr]$$
where $\bigl[\,_i\tilde{x},\,_i\tilde{y},\,_i\tilde{z}\bigr]$,
$i=1,\ldots,k$, are $k$ different elements of $W_{G'}$. By 
theorem \ref{theoabasis} and proposition \ref{LGbasis} we deduce that
$\mathbb{B}_{\phi}(\alpha)=0$ and as the elements of $W_{G'}$
form a basis the coefficient $a_i$,
$i=1,\ldots,k$, have to be zero. Thus there are no relations
among the elements of $N_5$ and between $N_5$ and $N_6\cup N_7$. \\
For the set $N_6$ we apply again lemma \ref{lemmalinind} with $M\cup
Q=N_6$ and $N=N_7$ which shows that there are no relations among the
elements of $N_6$ and between the elements of $N_6$ and $N_7$, besides
the two torsion.\\
Finally the set $N_7$ remains. There again lemma \ref{lemmalinind} can
be applied where $M=N_7$ and $Q=N=\emptyset$ to show that the elements
are linearly independent. 
\endBeweis
\begin{theo}\label{basisB_G}
The set $\cup_{i=1}^7N_i$ is a basis for $\mathbb{B}_G$.
\end{theo}
\beginBeweis
By proposition \ref{spanB_G} the set $\cup_{i=1}^7N_i$ generates the
group $\mathbb{B}_G$ and by proposition \ref{linunabB_G} there are no
relations among the elements
of $\cup_{i=1}^7N_i$, except for the two torsion and thus they form a basis.
\endBeweis
\begin{theo}\label{isoB_G}
The group of singular orbit data $\mathbb{B}_G$ is isomorphic to 
$$\mathbb{Z}^{|T^+_G|}\oplus (\mathbb{Z}/2\mathbb{Z})^{|S_G|-n_{K'}}.$$
\end{theo}
\beginBeweis
We will prove this theorem by introducing a one to one correspondence
between the elements of $T^+_G\cup S_G$, without the images under the map
$\hat{\psi}$ of the $n_{K'}$ generators of $K'$, and the
elements of $\cup_{i=1}^7N_i$.  
\begin{enumerate}
\item To every element $x\in[G,G]$ with $\hat{x}\in T^+_G\cup S_G$
corresponds the element $\bigl[\hat{x}\bigr]\in N_1$ and vice versa. 
The elements of $T^+_G$ give rise to copies of $\mathbb{Z}$ and the
elements of $S_G$ give rise to copies of $\mathbb{Z}/2\mathbb{Z}$.

\item A conjugacy class $\hat{x}_i\in T^+_G$, $i>0$, with
$\tilde{x}\neq\tilde{x}^{-1}$ corresponds to the element
$\bigl[\hat{x}_0,\hat{x}_i^{-1}\bigr]\in N_2$ and vice versa. The
elements of $N_2$ give all rise to copies of $\mathbb{Z}$.

\item Let $\hat{x}_i\in T^+_G$, $i>0$, with $\tilde{x}=\tilde{x}^{-1}$
and $x\not\in[G,G]$, then there are two possibilities.
\begin{enumerate}
\item $\hat{x}_0\in T^+_G$; Thus the conjugacy class $\hat{x}_i$ corresponds
to the element $\bigl[\hat{x}_0,\hat{x}_i\bigr]\in N_3$ and vice
versa. The elements of $N_3$ give rise to copies of $\mathbb{Z}$.

\item $\hat{x}_0\in S_G$; Thus the conjugacy class $\hat{x}_i$ corresponds
to the element $\bigl[\hat{x}_0,\hat{x}_i\bigr]\in N_4$. These
elements of $N_4$ give also rise to copies of $\mathbb{Z}$. For the
other elements of $N_4$ see the next item.
\end{enumerate}

\item Let $\hat{x}_i\in S_G$, $i>0$, then it follows $\hat{x}_0\in
S_G$ and the corresponding element $\bigl[\hat{x}_0,\hat{x}_i\bigr]\in
N_4$ gives rise to a copy of $\mathbb{Z}/2\mathbb{Z}$. With this we
found a correspondence with every element of $N_4$.

\item Now we want to look at the conjugacy classes $\hat{x}_0\in T^+_G$ with
$\tilde{x}\neq\tilde{x}^{-1}$. Between the sets $T^+_G$ and $T^+_{G'}$
we have the following relation:\\
$\#\bigl|\bigl\{\hat{x}_0\in T^+_G\,|\,\tilde{x}\neq
\tilde{x}^{-1}\bigr\}\bigr|=\#|T^+_{G'}|$. 

On the other hand by the construction of $N_5$ we have: \\
$\#|W_{G'}|=\#|N_5|$.

By theorem \ref{theoaisom} we then obtain the equality:\\
$\#\bigl|\bigl\{\hat{x}_0\in T^+_G\,|\,\tilde{x}\neq
\tilde{x}^{-1}\bigr\}\bigr|=\#|T^+_{G'}|=\#|W_{G'}|=\#|N_5|$.

The elements of $N_5$ give rise to copies of $\mathbb{Z}$ because
$\mathbb{B}_{\phi}$ is a homomorphism and the images have infinite
order.

Thus we have seen that all the elements $\hat{x}_0\in T^+_G$ with
$\tilde{x}\neq\tilde{x}^{-1}$ give rise to copies of $\mathbb{Z}$ in
$\mathbb{B}_G$. 

\item Recall the notation we introduced to define the sets $N_6$ and
$N_7$. In this notation an element $\tilde{x}\in H_2'$ can be written
as a product $$\tilde{x}=\,_{i_1}\tilde{z}\cdots\,_{i_t}\tilde{z}$$
with $i_j=1,\ldots,n_{G'}$, $1\leq j\leq t$, $t\geq1$, and we have
$\hat{x}_0=\mathbb{B}_{\hat{\psi}}(\,_{i_1}\tilde{z}\cdots\,_{i_t}\tilde{z})=
\,_{i_1}\hat{z}_0\cdots\,_{i_t}\hat{z}_0$. Note that the images under
the map $\hat{\psi}$ of the $n_{K'}$ generators of $K'$
are denoted by $_i\hat{z}_0$, $i=1,\ldots,n_{K'}$. 

Let now the conjugacy classes $\hat{x}_0\in S_G$, but  
$\hat{x}_0\neq\,_i\hat{z}_0$ with $i=1,\ldots,n_{K'}$ (i.e.,
$\hat{x}_0$ is not the image under the map $\hat{\psi}$
of one of the $n_{K'}$ generators of $K'$). Since 
$\tilde{x}$ is an element of $M_{G'}$ but not a generator of $K'$ we
have the following presentation for $\tilde{x}$.
$$\tilde{x}=\,_{i_1}\tilde{z}\cdots\,_{i_t}\tilde{z}$$
with $i_j=1,\ldots,n_{K'}$, $1\leq j\leq t$ and $t>1$. The element
which corresponds to 
$\hat{x}_0$ is now $\bigl[\,_{i_1}\hat{z}_0,\ldots,\,_{i_t}\hat{z}_0,
\hat{x}_0\bigr]\in N_6$ and this element has order two. The remaining
elements of $N_6$ are covered by the next item.

\item The last conjugacy classes which remain from the set $T^+_G\cup S_G$,
without the images under the map $\hat{\psi}$ of the
$n_{K'}$ generators of $K'$, are the classes $\hat{x}_0\in T^+_G$ with
$\tilde{x}=\tilde{x}^{-1}$ and $\tilde{x}\neq 1$ (i.e.,
$\tilde{x}\in H_2'-M_G'$). With the notation of the previous item we
have again two cases. 
\begin{enumerate}
\item $t>1$; $\hat{x}_0$ corresponds to the element
$\bigl[\,_{i_1}\hat{z}_0,\ldots,\,_{i_t}\hat{z}_0,\hat{x}_0\bigr]\in
N_6$ which gives rise to a copy of $\mathbb{Z}$. These are now all
elements of $N_6$. 

\item $t=1$; Then we have $\hat{x}_0=\,_{i_1}\hat{z}_0$ and thus 
$i_1=n_{K'}+1,\ldots,n_{G'}$ and the conjugacy class corresponds to
the element $\bigl[\hat{x}_0,\hat{x}_0\bigr]\in N_7$ and vice
versa. Note that the elements of $N_7$ have all infinite order.
\end{enumerate}
\end{enumerate}
With this we have shown that every element of $T^+_G\cup S_G$, without
the images under the map $\hat{\psi}$ of the $n_{K'}$
generators of $K'$, corresponds to exactly one element in
$\cup_{i=1}^7N_i$ with the appropriate order. 
\endBeweis
\begin{cor}\label{trivial}
The group $\mathbb{B}_G$ is trivial if and only $G\cong C_2$ or $G$ is
trivial.
\end{cor}
\beginBeweis
The case where $G$ is trivial is trivial.
If $G\cong C_2$ then by example~1, $\mathbb{B}_G$ is trivial.

On the other hand let $\mathbb{B}_G\cong 0$, then $T_G=\emptyset$ and
$|S_G|=n_{K'}$. From this we deduce that $G$ consists only of elements
which are conjugate to their inverse, i.e., $|\hat{G}|-1=|S_G|$ and thus
$n_{K'}=n_{G'}$. Moreover we have $|G'|=2^{n_{G'}}$. From this we
conclude:
$$2^{n_{G'}}-1=|G'|-1\leq|\hat{G}|-1=|S_G|=n_{G'}.$$
This equation yields two cases either $n_{G'}=0$ or $n_{G'}=1$. If
$n_{G'}$ is zero then $|S_G|=0$ and the group $G$ is trivial. If
$n_{G'}$ is one then $G/[G,G]\cong C_2$ and $|\hat{G}|=2$. Thus $[G,G]$
has to be trivial and $G\cong C_2$.
\endBeweis

\section{Relation with Cobordism}\label{cobordism}
Before we can talk about the relation with $G$-equivariant cobordism
we give its definition. 
\begin{defin}
Let $M_1$ (resp. $M_2$) be a compact, oriented, connected Riemann
surface with smooth $G$-action $\kappa_1:G\rightarrow\text{\it
Diffeo}_+(M_1)$ (resp. $\kappa_2:G\rightarrow\text{\it Diffeo}_+(M_2)$) 
We say that $\kappa_1$ is {\bf $G$-equivariant cobordant} to
$\kappa_2$, written 
$\kappa_1\sim\kappa_2$, if there exists a smooth, compact, oriented, connected
3-manifold $V$ and a smooth $G$-action $\Phi$ on $V$ such that
\begin{enumerate}
\item The boundary of $V$ is the disjoint union of $M_1$ and $-M_2$,
$\partial(V)=M_1\cup-M_2$. The notation $-M_2$ denotes $M_2$ with opposite
orientation. The orientations on $M_1$ and $-M_2$ coincide with the
one induced by $V$. 
\item $\Phi$ restricted to $\partial(V)$ agrees with $\kappa_1\cup\kappa_2$.
\end{enumerate}
We also say that $\kappa_1$ is zero $G$-equivariant cobordant, written
$\kappa_1\sim0$, if $\partial(V)=M_1$.
$\Omega_G$ will denote the set of $G$-equivariant cobordism classes
and a class will be denoted by $(\kappa,M)$.
\end{defin}
The set $\Omega_G$ forms an Abelian group where the addition is given
as for $\mathbb{B}_G$ in section \ref{sec:group} by the
$G$-equivariant connected sum. In 
\cite{Gr3} the author shows that two $G$-actions which are cobordant
have the same singular orbit data, thus there is a well defined
homomorphism $\chi:\Omega_G\rightarrow\mathbb{B}_G$ which sends every
class $(\kappa,M)$ to its singular orbit data. Moreover the map
$\chi$ is surjective. Indeed, take any $G$-action which represents a
given singular orbit data $\alpha\in\mathbb{B}_G$, the corresponding
cobordism class will then be 
mapped to the same singular orbit data $\alpha$ by $\chi$. In the same paper it
is also shown that the kernel of $\chi$ consists only of the
cobordism classes of free $G$-actions. The next proposition proves
that the subgroup of free actions is isomorphic to
$H_2(G;\mathbb{Z})$. The proof is an easy consequence of a spectral
sequence described in Conner and Floyd's book \cite{CoFl}.
\begin{prop}\label{kernel}
The kernel of the map $\chi$ is isomorphic to $H_2(G;\mathbb{Z})$.
\end{prop}
\beginBeweis
In \cite[Theorem 20.4]{CoFl} Conner and Floyd prove that the subgroup
of cobordism classes of free $G$-actions is isomorphic to $MSO_2(BG)$,
the cobordism homology of the classifying space of $G$. 

On the other hand there is a spectral sequence $\{E^r_{p,q}\}$ with
$E^2_{p,q}=H_p(G;MSO_q)$ and whose $E^{\infty}$-term is associated to
a filtration of $MSO_*(BG)$. 
It turns out that for $MSO_2(BG)$ the $E^2$-term already stabilizes and
as $MSO_0\cong\mathbb{Z}$ and $MSO_1\cong MSO_2\cong0$ we have
$$MSO_2(BG)\cong\sum_{p+q=2}H_p(G;MSO_q)\cong H_2(G;\mathbb{Z}).$$
\endBeweis
We can now conclude that the $G$-equivariant cobordism group
$\Omega_G$ of surface diffeomorphisms is an Abelian group extension of
$\mathbb{B}_G$ by $H_2(G;\mathbb{Z})$.
\begin{cor}
Every $G$-action is cobordant to a free action if and only if $G\cong
C_2$ or $G$ is the trivial group. 
\end{cor}
\beginBeweis
By corollary \ref{trivial}, $\mathbb{B}_G$ is the trivial group if and
only if $G\cong C_2$ or $G$ is trivial. Thus the corollary follows as
the kernel of $\chi$ consists of the cobordism classes of the free
$G$-actions. 
\endBeweis
\begin{cor}
If there is no torsion in the group $\mathbb{B}_G$, then
$\Omega_G$ is isomorphic to the direct sum of $\mathbb{B}_G$ with
$H_2(G;\mathbb{Z})$. 
\end{cor}
\beginBeweis
The group $\mathbb{B}_G$ consists only of copies of $\mathbb{Z}$ and
$\Omega_G$ surjects onto this group. Thus as $\Omega_G$ is Abelian and
finitely generated the short exact sequence
$0\rightarrow
H_2(G;\mathbb{Z})\rightarrow\Omega_G\rightarrow\mathbb{B}_G\rightarrow
0$
splits.
\endBeweis

\section{Relation to Representation Theory of finite
Groups}\label{reptheory}
In section \ref{sec:group} we introduced the $G$-signature of Atiyah
and Singer $\eta:\mathbb{B}_G\rightarrow R_{\mathbb{C}}(G)$ (equation
(\ref{eta}))  and the
$G$-signature $\theta:\mathbb{B}_G\rightarrow R_{\mathbb{C}}(G)/E_G$,
but we referred to this section for a short discussion of the subgroup
$E_G$ and the properties of~$\theta$. 

The subgroup $E_G$ is defined as follows:
$$E_G := \bigl< Ind^G_H\rho_0\,|\,H\leq G \bigr>.$$
The representation $\rho_0$ denotes always the one dimensional trivial
representation. The map $\theta$ is then just the map $\varphi$
followed by the surjection on to the quotient and it turns out that
the resulting map is a homomorphism (see \cite[section 4]{Gr3}). 

It is proven in \cite[theorem 21]{Gr3} that the map $\theta$ is injective
on the copies of $\mathbb{Z}$ in $\mathbb{B}_G$ and moreover
that $\theta(\beta)\neq\overline{\theta(\beta)}$ for every element of
infinite order $\beta\in\mathbb{B}_G$ and
$\theta(\alpha)=\overline{\theta(\alpha)}$ for every element of order
two $\alpha\in\mathbb{B}_G$ \cite[proposition 24]{Gr3}. By
the way $\eta$ is defined in equation (\ref{eta}) it follows that
$\eta$ is also injective on the elements of infinite order but zero on
the elements of order two. 

What can we say about $\theta$? Is it also zero on the elements of
order two? The subgroup $E_G<R_{\mathbb{C}}(G)$ is contained in
$R_{\mathbb{Q}}(G)$ which shows that $\theta(\alpha)$ might be zero
for elements $\alpha$ of order two as
$\theta(\alpha)=\overline{\theta(\alpha)}$ in this case.

For $G=S_3$ (see example 5), the symmetric group on three
letters, we have $\mathbb{B}_{S_3}\cong C_2$ and on the other hand
$E_G=R_{\mathbb{Q}}(G)=R_{\mathbb{C}}(G)$. From the second fact we
deduce that the map $\theta$ is the zero map. This phenomenon where
for a finite group $G$ its {\it geometric
representation theory} $\mathbb{B}_G$ and its
representation ring $R_{\mathbb{C}}(G)$ both have special properties
fits into a broader picture. 
\begin{prop}
The group $\mathbb{B}_G$ contains only elements of order two if and
only if every complex character of $G$ has values in $\mathbb{R}$.
\end{prop}
\beginBeweis
The group $\mathbb{B}_G$ consists only of elements of order two if and
only if every element of $G$ is conjugate to its inverse. 

If every element $a$ of $G$ is conjugate to its inverse, then by the
formula $\overline{\chi(a)}=\chi(a^{-1})$, $\chi$ a character, every
character has values in $\mathbb{R}$. On the other hand if every
character $\chi$ has values in $\mathbb{R}$, then the characters have
the same 
value for an element $a$ of $G$ and their inverse $a^{-1}$. But the
characters form a basis for the vector space of class functions and as
such have to separate conjugacy classes. Thus every element has to be
conjugate to its inverse. (See also Serre's book \cite{Se}.)
\endBeweis
The proposition doesn't say anything about $E_G$, so we don't know in
general whether the map $\theta$ is zero or not in this case.

The advantage of our approach is that in the case where $\theta$ is
not injective,
we can try to find a smaller subgroup $E_G'<E_G$ such that the new map
$$\theta':\mathbb{B}_G\rightarrow R_{\mathbb{C}}(G)/E_G'$$ 
is still a homomorphism and in addition becomes injective. We want to
illustrate this idea with the example $G=S_3$.

Let $\chi_0$ be the trivial representation, $\chi_1$ the one
dimensional non-trivial representation and $\chi_2$ the
two dimensional irreducible representation of $S_3$. Recall the notation from
example 5 and let 
$$E'_{S_3}=\bigl<Ind^{S_3}_{<a>}\rho_0\,,\,Ind^{S_3}_{<1>}\rho_0\,,\,Ind^{S_3}_{S_3}\rho_0\bigr>
=\bigl<\chi_0+\chi_1\,,\,\chi_0+\chi_1+2\chi_2\,,\,\chi_0\bigr>.$$
Then we have that the class of $\chi_2$ generates
$R_{\mathbb{C}}S_3/E'_{S_3}\cong C_2$. We know
that $\mathbb{B}_{S_3}\cong\bigl<\bigl[\hat{a}\bigr]\bigr>\cong C_2$; thus it
remains to show that $\theta':\mathbb{B}_{S_3}\rightarrow
R_{\mathbb{C}}S_3/E'_{S_3}$ is non zero on $\bigl[\hat{a}\bigr]$.

Let $\phi_{[\hat{a}]_{S_3}}$ be an embedding of $S_3$ into $\Gamma_3$
with singular 
orbit data $[\hat{a}]_{S_3}$. Then by the method of \cite[proposition
17]{Gr3} we find that $\varphi(\phi_{[\hat{a}]_{S_3}})=\chi_0+\chi_2$
and thus 
$\theta'$ maps $\bigl[\hat{a}\bigr]$ to the class of $\chi_2$ and
$\theta'$ is an isomorphism.

Note that $E'_{S_3}$ consists of the minimal relations in order to make
$\theta'$ a group homomorphism. It consists of the representations
coming from the free actions and the cancelling pair
$\bigl[\hat{a},\hat{a}\bigr]_{S_3}$.

In this situation the maps
$\mathbb{B}_i:\mathbb{B}_{C_3}\rightarrow\mathbb{B}_{S_3}$,
$\widetilde{Ind}_{C_3}^{S_3}:\widetilde{R_{\mathbb{C}}C_3}\rightarrow
R_{\mathbb{C}}S_3/E'_{S_3}$ and $\theta'$ still commute. This follows
from the fact that  
$$Ind_{C_3}^{S_3}(\varphi(\phi_{[a,a,a]_{C_3}}))\equiv
Ind_{C_3}^{S_3}\rho_1=\chi_2\equiv\varphi(\phi_{[\hat{a}]_{S_3}})
\equiv\varphi(\phi_{\mathbb{B}_i([a,a,a]_{C_3})})$$ 
where the equivalence is taken modulo $E'_{S_3}$ and $\rho_1$ denotes a
one dimensional faithful representation of $C_3$. With $\phi_{[a,a,a]_{C_3}}$
and $\phi_{\mathbb{B}_i([a,a,a]_{C_3})}$ we denote embeddings of
$C_3$ and $S_3$ respectively, into
some mapping class groups with singular orbit data $[a,a,a]_{C_3}$ and
$\mathbb{B}_i([a,a,a]_{C_3})$ respectively.

\end{document}